\input amstex
\documentstyle{amsppt}
%\magnification 1200
\NoBlackBoxes

\TagsOnRight

\def\cal{\Cal}
\def\AA{{\cal A}}
\def\BB{{\cal B}}

\def\HH{{\cal H}}
\def\MM{{\cal M}}

\def\Z{{\Bbb Z}}

\def\R{{\Bbb R}}

\def\Q{{\Bbb Q}}
\def\e{{\epsilon}}

\def\n{\noindent}
\def\part{{\partial}}
\def\dudtau{{\part u\over \part \tau}}
\def\dudt{{\part u\over \part t}}

\rightheadtext{Spectral invariants} \leftheadtext{
Yong-Geun Oh }

\topmatter
\title
Construction of spectral invariants \\
of Hamiltonian paths on \\
closed symplectic manifolds
\endtitle
\author
Yong-Geun Oh\footnote{Partially supported by the NSF Grant \#
DMS-9971446, by \# DMS-9729992 in the Institute for Advanced
Study, Vilas Associate Award in the University of Wisconsin
and by a grant of the Korean Young Scientist Prize
\hskip8.5cm\hfill}
\endauthor
\address
Department of Mathematics, University of Wisconsin, Madison, WI
53706, ~USA \& Korea Institute for Advanced Study, Seoul, Korea;
oh\@math.wisc.edu
\endaddress

\abstract  In this paper, we develop a mini-max theory of the
action functional over the semi-infinite cycles via the chain
level Floer homology theory and construct spectral invariants of
Hamiltonian diffeomorphisms on arbitrary, especially on {\it
non-exact and non-rational}, compact symplectic manifold
$(M,\omega)$. To each given time dependent Hamiltonian function
$H$ and quantum cohomology class $ 0 \neq a \in QH^*(M)$, we
associate an invariant $\rho(H;a)$ which varies continuously over
$H$ in the $C^0$-topology. This is obtained as the mini-max value
over the semi-infinite cycles whose homology class is `dual' to
the given quantum cohomology class $a$ on the covering space
$\widetilde \Omega_0(M)$ of the contractible loop space
$\Omega_0(M)$. We call them the {\it Novikov Floer cycles}. We
apply the spectral invariants to the study of Hamiltonian
diffeomorphisms in sequels of this paper.

We assume that $(M,\omega)$ is strongly semi-positive in this paper,
which will be removed in a sequel to this paper.

\endabstract


\endtopmatter
\document

\centerline {\it Dedicated to Alan Weinstein in honor of his 60-th
birthday}

\bigskip

\centerline{\bf Contents} \medskip

\n \S1. Introduction and the main results
\smallskip

\n \S2. The action functional and the action spectrum
\smallskip

\n \S3. Quantum cohomology in the chain level
\smallskip

\n \S4. Filtered Floer homology  \par
\smallskip

\n \S5. Construction of the spectral invariants of Hamiltonians
\smallskip
\quad 5.1. Overview of the construction \par
\quad 5.2. Finiteness; the linking property of semi-infinite cycles  \par

\smallskip

\n \S6. Basic properties of the spectral invariants
\smallskip
\quad  6.1. Proof of the symplectic invariance \par
\smallskip
\quad  6.2. Proof of the triangle inequality
\smallskip

\n \S7. The rational case: proof of the spectrality
\par
\smallskip

\n \S8. Remarks on the transversality
\smallskip

\n {\it Appendix}: Continuous quantum cohomology
\smallskip

\vskip0.5in

\head
{\bf \S 1. Introduction and the main results}
\endhead

The group $\HH am(M,\omega)$ of (compactly supported) Hamiltonian
diffeomorphisms of the symplectic manifold $(M,\omega)$ carries
a remarkable invariant norm  defined by
$$
\aligned
\|\phi\| & = \inf_{H \mapsto \phi} \|H\| \\
\|H\| & = \int_0^1(\max H_t - \min H_t)\, dt
\endaligned
\tag 1.1
$$
which was introduced by Hofer [Ho].
Here $H\mapsto \phi$ means that $\phi$ is the time-one map $
\phi_H^1$ of the Hamilton's equation $\dot x =X_H(x)$ of
the Hamiltonian $H: [0,1] \times M \to \R$, where the Hamiltonian vector
field is defined by
$$
\omega(X_H, \cdot) = dH. \tag 1.2
$$
This norm can be easily defined on arbitrary symplectic manifolds
although proving non-degeneracy is a  non-trivial matter (See
[Ho], [Po1] and [LM] for its proof of increasing generality. See
also [Ch] for a Floer theoretic proof and [Oh3] for a simple proof
of the non-degeneracy in tame symplectic manifolds).

On the other hand Viterbo [V] defined another invariant norm on
$\R^{2n}$. This was defined by considering the graph of the
Hamiltonian diffeomorphism $\phi: \R^{2n} \to \R^{2n}$ and
compactifying the graph in the diagonal direction in $\R^{4n} =
\R^{2n}\times \R^{2n}$ into $T^*S^{2n}$.  He then applied the
critical point theory of generating functions of the Lagrangian
submanifold $\operatorname{graph}\phi \subset T^*S^{2n}$ which he
developed on the cotangent bundle $T^*N$ of the arbitrary compact
manifold $N$. To each cohomology class $a \in H^*(N)$, Viterbo
associated certain homologically essential critical values of
generating functions of any Lagrangian submanifold $L$ Hamiltonian
isotopic to the zero section of $T^*N$ and proved that they depend
only on the Lagrangian submanifold but not on the generating
functions, at least up to normalization.

The present author [Oh1,2] and Milinkovi\'c [MO1,2, M] developed a
Floer theoretic approach to construction of Viterbo's invariants
using the canonically defined action functional on the space of
paths, utilizing the observation made by Weinstein [W] that the
action functional is a generating function of the given Lagrangian
submanifold defined on the path space. This approach is canonical
{\it including normalization} and provides a direct link between
Hofer's geometry and Viterbo's invariants in a transparent way.
One of the key points in our construction in [Oh2] is the emphasis
on the usage of the existing group structure on the space of
Hamiltonians defined by
$$
(H,K) \mapsto H\# K: = H + K \circ (\phi_H^t)^{-1} \tag 1.3
$$
in relation to the pants product and the triangle inequality.
However we failed to fully exploit this structure and fell short
of proving the triangle inequality at the time of writing [Oh1,2].

This construction can be carried out for the Hamiltonian
diffeomorphisms as long as the action functional is single valued,
e.g., on {\it weakly-exact} symplectic manifolds. Schwartz [Sc]
carried out this construction in the case of {\it symplectically
aspherical} $(M,\omega)$, i.e., for $(M,\omega)$ with
$c_1|_{\pi_2(M)} = \omega|_{\pi_2(M)} = 0$. Among other things he
proved the triangle inequality for the invariants constructed
using the notion of Hamiltonian fibration and (flat) symplectic
connection on it. It turns out that the proof of this triangle
inequality [Sc] is closely related to the notion of the $K$-area
of the Hamiltonian fibration [Po2] with connections [GLS], [Po2],
especially to the one with fixed monodromy studied by Entov [En1].
In this context, the choice of the triple $(H,K; H\#K)$ we made in
[Oh2] can be interpreted as the one which makes infinity the
$K$-area of the corresponding Hamiltonian fibration over the
Riemann surface of genus zero with three punctures equipped with
the given monodromy around the punctures. Entov [En1] develops a
general framework of Hamiltonian connections with fixed boundary
monodromy and relates the $K$-area with various quantities of the
given monodromy which are of the Hofer length type. This framework
turns out to be particularly useful for our construction of
spectral invariants in the present paper.

On non-exact symplectic manifolds, the action functional is not
single valued and the Floer homology theory has been developed as
a circle-valued Morse theory or a Morse theory on a covering space
$\widetilde \Omega_0(M)$ of the space $\Omega_0(M)$ of
contractible (free) loops on $M$ in the literature related to
Arnold's conjecture which was initiated by Floer himself [Fl]. The
Floer theory now involves quantum effects and uses the Novikov
ring in an essential way [HoS]. The presence of quantum effects
and {\it denseness of the action spectrum} in $\R$ (as in
non-rational symplectic manifolds), had been the most serious
obstacle that plagued the study of {\it family of Hamiltonian
diffeomorphisms}, until the author [Oh4] developed a general
framework of the mini-max theory over natural semi-infinite cycles
on the covering space $\widetilde \Omega_0(M)$ which we call the
{\it Novikov Floer cycles}. In the present paper, we will exploit the
`finiteness' condition in the definitions of the Novikov ring and
the Novikov Floer cycles in a crucial way for the proofs of various
existence results of pseudo-holomorphic curves that are needed in
the proofs of the axioms of spectral invariants and nondegeneracy
of the norm that we construct [Oh8]. Although the Novikov ring is
essential in the definition of the Floer homology and the quantum
cohomology in the literature, as far as we know, it is the first
time for the finiteness condition to be explicitly used beyond the
purpose of giving the definition of the quantum cohomology and the
Floer homology.

A brief description of the setting of the Floer theory [HoS] is
in order, partly to fix our convention:
Let $(\gamma,w)$ be a pair of
$\gamma \in \Omega_0(M)$ and $w$ be a disc bounding $\gamma$. We
say that $(\gamma,w)$ is {\it $\Gamma$-equivalent} to
$(\gamma,w^\prime)$ iff
$$
\omega([w'\# \overline w]) = 0 \quad \text{and }\, c_1([w'\#
\overline w]) = 0 \tag 1.4
$$
where $\overline w$ is the map with opposite orientation on the
domain and $w'\# \overline w$ is the obvious glued sphere. Here
$\Gamma$ stands for the group
$$
\Gamma = {\pi_2(M)\over \text{ker\ } (\omega|_{\pi_2(M)}) \cap
\text{ker\ } (c_1|_{\pi_2(M)})}.
$$

We denote by $[\gamma,w]$ the
$\Gamma$-equivalence class of $(\gamma,w)$ and by
$\widetilde \Omega_0(M)$ the set of $\Gamma$-equivalence classes.
Let $\pi: \widetilde
\Omega_0(M) \to \Omega_0(M)$ the canonical projection. We
call $\widetilde \Omega_0(M)$ the {\it $\Gamma$-covering space} of
$\Omega_0(M)$. The action functional $\AA_0: \widetilde
\Omega_0(M) \to \R$ is defined by
$$
\AA_0([\gamma,w]) = -\int w^*\omega. \tag 1.5
$$
Two $\Gamma$-equivalent pairs $(\gamma,w)$ and $(\gamma,w^\prime)$
have the same action and so the action is well-defined on
$\widetilde\Omega_0(M)$. When a one-periodic Hamiltonian
$H:(\R/\Z) \times M \to \R$ is given, we consider the functional $\AA_H:
\widetilde \Omega(M) \to \R$ by
$$
\AA_H([\gamma,w])= -\int w^*\omega - \int H(t, \gamma(t))dt.
\tag 1.6
$$
Our convention is chosen to be consistent with the classical mechanics
Lagrangian on the cotangent bundle with the symplectic form
$$
\omega_0 = - d\theta, \quad \theta = \sum_i p_i dq^i
$$
when (1.2) is adopted as the definition of Hamiltonian vector
field. See the remark in the end of this introduction on other
conventions in the symplectic geometry. The conventions in the
present paper coincide with our previous papers [Oh1,2,4] and
Entov's [En1,2] but different from many other literature on the
Floer homology one way or the other. (There was a sign error in
[Oh1,2] when we compare the Floer complex and the Morse complex
for a small Morse function, which was rectified in [Oh4]. In our
convention, the  positive gradient flow of $\e f$ corresponds to
the negative gradient flow of $\AA_{\e f}$.)

The mini-max theory of this action functional on the
$\Gamma$-covering space has been implicitly used in the proof of
Arnold's conjecture. Recently the present author has further
developed this mini-max theory via the Floer homology and applied
it to the study of Hofer's geometry of Hamiltonian diffeomorphism
groups [Oh4]. We also outlined construction of spectral invariants
of Hamiltonian diffeomorphisms of the type [V], [Oh2], [Sc] on
arbitrary non-exact symplectic manifolds for the {\it classical}
cohomological classes. The main purpose of the present paper is to
further develop the chain level Floer theory introduced in [Oh4]
and to carry out  construction of spectral invariants for
arbitrary quantum cohomology classes.  The organization of the
paper is now in order.

In \S 2, we briefly review various facts related to the action
functional and its action spectrum. Some of these may be known to
the experts, but precise details for the action functional on the
covering space $\widetilde \Omega_0(M)$ of general $(M,\omega)$
first appeared in our paper [Oh5] especially concerning the
normalization and the loop effect on the {\it action spectrum}: We
define the action spectrum of $H$ by
$$
\text{Spec}(H): = \{ \AA_H([z,w]) \in \R \mid [z,w] \in \widetilde
\Omega_0(M), \, d\AA_H([z,w]) = 0\}
$$
i.e., the set of critical values of $\AA_H: \widetilde \Omega_0(M)
\to \R$. In [Oh5], we have shown that once we normalize the
Hamiltonian $H$ on compact $M$ by
$$
\int_M H_t\, d\mu = 0
$$
with $d\mu$ the Liouville measure, $\text{Spec}(H)$ depends only
on the equivalence class $\widetilde \phi = [\phi,H]$ (see \S 2 for
the definition) and so $\text{Spec}(\widetilde \phi) \subset \R$ is
a well-defined subset of $\R$ for each $\widetilde \phi \in
\widetilde{\HH am}(M,\omega)$. Here
$$
\pi: \widetilde{\HH am}(M,\omega) \to \HH am(M,\omega)
$$
is the universal covering space of $\HH am(M,\omega)$.
This kind of
normalization of the action spectrum is a crucial point for
systematic study of the spectral invariants of the Viterbo type in
general. Schwarz [Sc] previously proved that in the aspherical
case where the action functional is single valued already on
$\Omega_0(M)$, this normalization can be made on $\HH
am(M,\omega)$, not just on $\widetilde{\HH am}(M,\omega)$.

In \S 3, we review the quantum cohomology and its Morse theory
realization of the corresponding complex. We emphasize the role of
the Novikov ring in relating the quantum cohomology and the Floer
homology and the reversal of upward and downward Novikov rings in
this relation. In \S 4, we review the standard operators in the
Floer homology theory and explain the filtration naturally present
in the Floer complex and how it changes under the Floer chain map.
In \S 5,  we give the definition of our spectral invariants for
the Hamiltonian functions $H$, and prove finiteness of the
mini-max values $\rho(H;a)$. In \S 6, we prove all the basic
properties of the spectral invariants. We summarize these into the
following theorem. We denote by $C_m^\infty([0,1]\times M)$ the
set of normalized continuous functions on $[0,1] \times M$.
Noting that there is a one-one correspondence between the set
$C_m^\infty([0,1]\times M)$ and the set of {\it Hamiltonian paths}
$$
\lambda= \phi_H: t \in [0,1] \mapsto \phi_H^t \in \HH
am(M,\omega),
$$
one may equally consider $\rho(H;a)$ as an invariant attached to the
Hamiltonian path $\phi_H$.

\proclaim{Theorem I} Let $(M,\omega)$ be arbitrary closed
symplectic manifold. For any given quantum cohomology class $0
\neq a \in QH^*(M)$, we have a continuous function denoted by
$$
\rho =\rho(H; a): C_m^\infty([0,1] \times M) \times QH^*(M) \to \R
$$
such that they satisfy the following axioms: Let $H, \, F  \in
C_m^\infty([0,1] \times M)$ be smooth Hamiltonian
functions and $a \neq 0 \in QH^*(M)$. Then $\rho$ satisfies the
following axioms:

\roster \item {\bf (Projective invariance)} $\rho(H;\lambda a) =
\rho(H;a)$ for any $0 \neq \lambda \in \Q$. \item {\bf
(Normalization)} For $a = \sum_{A \in \Gamma} a_A q^{-A} $, we
have $\rho(\underline 0;a) = v(a)$ where $\underline 0$ is the
zero function and
$$
v(a): = \min \{\omega(-A) ~|~  a_A \neq 0 \} = - \max \{\omega(A)
\mid a_A \neq 0 \}. \tag 1.7
$$
is the (upward) valuation of $a$. \item {\bf (Symplectic
invariance)} $\rho(\eta^*H ;a) = \rho(H ;a)$ for any symplectic
diffeomorphism $\eta$ \item {\bf (Triangle inequality)} $\rho(H \#
F; a\cdot b) \leq \rho(H;a) + \rho(F;b) $ \item {\bf
($C^0$-continuity)} $|\rho(H;a) - \rho(F;a)| \leq \|H \# \overline
F\| = \|H - F\|$ where $\| \cdot \|$ is the Hofer's pseudo-norm on
$C_m^\infty([0,1] \times M)$. In particular, the function $\rho_a:
H \mapsto \rho(H;a)$ is $C^0$-continuous.
\endroster
\endproclaim

We will call the set
$$
\text{spec}(H): = \{ \rho(H;a) \mid a \in QH^*(M) \} \tag 1.8
$$
the {\it essential spectrum} of $H$.

Most of the properties stated in this theorem are direct analogs to the ones in
[Oh1,2] and [Sc]. {\it Except for the proof of finiteness of $\rho(H;a)$},
proofs of all of the properties are refinements of the arguments used
in [Oh2,4], [Sc].  In addition,  the proof of the triangle inequality uses
the concept of Hamiltonian fibration with fixed monodromy and the
$K$-area [Po2], [En1], which is an enhancement of the arguments
used in [Oh2], [Sc].

In the classical mini-max theory for the {\it indefinite}
functionals [Ra], [BnR], there was implicitly used the notion of
`semi-infinite cycles'
to carry out the mini-max procedure. There are two essential ingredients
needed to prove existence of actual critical values out of the
mini-max values: one is the finiteness
of the mini-max value, or the {\it linking property} of
the (semi-infinite) cycles associated to the class $a$ and the other is
to prove that the corresponding mini-max value
is indeed a critical value of the action functional.
{\it When the global gradient flow of the action functional
exists} as in the classical critical point theory [BnR],
this point is closely related to the
well-known Palais-Smale condition and the deformation lemma
which are essential ingredients needed to
prove the criticality of the mini-max value.
Partly because we do not have the global
flow, we need to geometrize all these classical mini-max procedures.
It turns out that the Floer homology theory in the chain level
is the right framework for this purpose.

In section 7, we will restrict to the {\it rational} case and
prove the following additional property of spectral invariants, {\it the spectrality axiom}.
We will study the non-rational cases elsewhere for which we expect the
same property holds, at least for the nondegenerate Hamiltonian functions, but
its proof seems to be much more nontrivial.

We now recall the definition of rational symplectic manifolds: Denote
$$
\Gamma_\omega := \{ \omega(A)~|~ A \in \pi_2(M)\} = \omega(\Gamma)
\subset \R
$$
and
$$
\hbox{\rm Spec}(H)= \cup_{z \in \text {Per}(H)}\text {Spec}(H;z).
$$
Recall that $\Gamma_\omega$ is either a discrete or a countable
dense subset of $\R$.

\definition{Definition 1.1} A symplectic manifold $(M,\omega)$ is
called {\it rational} if $\Gamma_\omega$ is discrete.
\enddefinition

\proclaim{Theorem II} {\bf (Spectrality Axiom)} Suppose that $(M,\omega)$ be rational.
Then
$\rho$ satisfies the following additional properties:

\roster \item For any smooth one-periodic Hamiltonian function $H:
S^1 \times M \to \R$, we have
$$
\rho(H;a) \in \text{\rm Spec}(H)
$$
for each given quantum cohomology class $0 \neq a \in QH^*(M)$.
\item For two smooth functions $H \sim K$ we have
$$
\rho(H;a) = \rho(K;a) \tag 1.9
$$
for all $a \in QH^*(M)$.
\endroster
\endproclaim

In particular, $\rho$ can be pushed down to the `universal covering space'
$\widetilde{\HH am}(M,\omega)$ of $\HH am(M,\omega)$
by putting $\rho(\widetilde \phi;a)$ to be the this
common value for $\widetilde \phi = [H]$. We call the subset
$\text{\rm spec}(\widetilde\phi)\subset \text{\rm Spec}(\widetilde
\phi)$ defined by
$$
\text{\rm spec}(\widetilde\phi) = \{ \rho(\widetilde\phi;a) \mid a
\in QH^*(M)\}
$$
the (homologically) {\it essential spectrum} of $\widetilde \phi$.
Then we have the following refined version of Theorem II for the
rational cases.

\proclaim{Theorem III} Let $(M,\omega)$ be rational and define the
map
$$
\rho: \widetilde{Ham}(M,\omega) \times QH^*(M) \to \R
$$
by $\rho(\widetilde\phi;a): = \rho(H;a)$. Let $\widetilde \phi, \,
\widetilde \psi \in \widetilde{Ham}(M,\omega)$ and $a \neq 0 \in
QH^*(M)$. Then $\rho$ satisfies the following axioms:

\roster \item {\bf (Spectrality)} For each $a\in QH^*(M)$,
$\rho(\widetilde \phi;a) \in \text{Spec}(\widetilde \phi)$. \item
{\bf (Projective invariance)} $\rho(\widetilde\phi;\lambda a) =
\rho(\widetilde\phi;a)$ for any $0 \neq \lambda \in \Q$. \item
{\bf (Normalization)} For $a = \sum_{A \in \Gamma} a_A q^{-A} $,
we have $\rho(\underline 0;a) = v(a)$ where $\underline 0$ is the
identity in $\widetilde{\HH am}(M,\omega)$ and
$$
v(a): = \min \{\omega(-A) ~|~  a_A \neq 0 \} = - \max \{\omega(A)
\mid a_A \neq 0 \}. \tag 1.10
$$
is the (upward) valuation of $a$. \item {\bf (Symplectic
invariance)} $\rho(\eta \widetilde \phi \eta^{-1};a) =
\rho(\widetilde \phi;a)$ for any symplectic diffeomorphism $\eta$
\item {\bf (Triangle inequality)} $\rho(\widetilde \phi \cdot
\widetilde \psi; a\cdot b) \leq \rho(\widetilde \phi;a) +
\rho(\widetilde \psi;b) $ \item {\bf ($C^0$-continuity)}
$|\rho(\widetilde \phi;a) - \rho(\widetilde \psi;a)| \leq
\|\widetilde \phi \circ \widetilde\psi^{-1} \| $ where $\| \cdot
\|$ is the Hofer's pseudo-norm on $\widetilde{Ham}(M,\omega)$. In
particular, the function $\rho_a: \widetilde \phi \mapsto
\rho(\widetilde \phi;a)$ is $C^0$-continuous. \item {\bf
(Monodromy shift)} Let $[h,\widetilde h] \in \pi_0(\widetilde G)$
act on $\widetilde{\HH am}(M,\omega) \times QH^*(M)$ by the map
$$
(\widetilde \phi,a) \mapsto (h\cdot \widetilde \phi, \widetilde
h^*a)
$$
where $\widetilde h^*a$ is the image of the (adjoint) Seidel's
action [Se] by $[h,\widetilde h]$ on the quantum cohomology
$QH^*(M)$. Then we have
$$
\rho([h,\widetilde h]\cdot (\widetilde \phi,a)) = \rho(\widetilde
\phi ;a) + I_\omega([h,\widetilde h]) \tag 1.11
$$
\endroster
\endproclaim

It would be an interesting question to ask whether these axioms
 characterize the spectral invariants
$\rho$. It is related to the question whether the graph of the
sections
$$
\rho_a: \widetilde \phi \mapsto \rho(\widetilde \phi;a); \quad
\widetilde{\HH am}(M,\omega) \to \frak{Spec}(M,\omega)
$$
can be split into other `branch' in a way that the other branch
can also satisfy all the above axioms or not. Here the {\it action
spectrum bundle} $\frak{Spec}(M,\omega)$ is defined by
$$
\frak{Spec}(M,\omega):= \bigcup_{\widetilde \phi\in \widetilde{\HH
am}(M,\omega)} \text{Spec}(\widetilde \phi) \subset \widetilde{\HH
am}(M,\omega) \times \R .
$$
We will investigate this question elsewhere.

 To get the main
stream of ideas in this paper without getting bogged down with
technicalities related with transversality question of various
moduli spaces, we  assume in this paper that $(M,\omega)$ is
strongly semi-positive in the sense of [Se], [En1]: A closed symplectic
manifold is called {\it strongly semi-positive} if there is no
spherical homology class $A \in \pi_2(M)$ such that
$$
\omega(A) > 0, \quad 2-n \leq c_1(A) \leq 0.
$$
Under this condition, the transversality problem concerning various moduli
spaces of pseudo-holomorphic curves is standard. We will not
mention this generic transversality question at all in the main
body of the paper unless it is absolutely necessary. In \S 7, we
will briefly explain how this general framework can be
incorporated in our proofs in the context of Kuranishi structure
[FOn] all at once. In Appendix, we introduce the notion of {\it
continuous quantum cohomology} and explain how to extend our
definition of spectral invariants to the continuous quantum
cohomology classes.

The present work is originated from a part of our paper entitled
``Mini-max theory, spectral invariants and geometry of the
Hamiltonian diffeomorphism group'' [Oh6] that has been circulated
since July, 2002. We isolate and streamline the construction part
of spectral invariants from [Oh6] in the present paper with some
minor corrections and addition of more details.
In particular, we considerably simplify the
definition of $\rho(H;a)$ from [Oh6] here. We leave the
application part of [Oh6] to a separate paper [Oh8] in which we
construct the homological norm of Hamiltonian diffeomorphism and
apply them to the study of geometry of Hamiltonian diffeomorphisms
on general compact symplectic manifolds.

Another application of the spectral invariants to the study of
length minimizing property of Hamiltonian paths is given by the
author [Oh7,8]. See also [En2], [EnP] for other interesting
applications of spectral invariants. In another sequel to this
paper, we will provide a description of spectral invariants in
terms of the Hamiltonian fibration.

We would like to thank the Institute for Advanced Study in
Princeton for the excellent environment and hospitality during our
participation of the year 2001-2002 program ``Symplectic Geometry
and Holomorphic Curves''. Much of the present work was finished
during our stay in IAS. We thank D. McDuff for some useful
communications in IAS. The final writing has been carried out
while we are visiting the Korea Institute for Advanced Study in
Seoul. We thank KIAS for providing an uninterrupted quiet time for
writing and excellent atmosphere of research.

We thank M. Entov and L. Polterovich for enlightening
discussions on spectral invariants and for explaining their
applications [En2], [EnP] of the spectral invariants
to the study of  Hamiltonian diffeomorphism group, and
Y. Ostrover for explaining his work from [Os] to us
during our visit of Tel-Aviv University. We also thank P.
Biran and L. Polterovich for their invitation to Tel-Aviv University
and hospitality.

\medskip

\n{\bf Convention.} \roster
\item The Hamiltonian vector field $X_f$ associated to a function
$f$ on $(M,\omega)$ is defined by $df = \omega(X_f,\cdot)$.
\item
The addition $F\# K$ and the inverse $\overline K$ on the set of
time periodic Hamiltonians $C^\infty(M \times S^1)$ are defined by
$$
\align
F\# G(x,t) & = F(x,t) + G((\phi_F^t)^{-1}(x),t) \\
\overline G(x,t) & = - G(\phi_G^t(x),t).
\endalign
$$
\endroster

There is another set of conventions which are used in the literature
(e.g., in [Po3]):
\roster
\item $X_f$ is defined by $\omega(X_f,\cdot) = -df$
\item The action functional has the form
$$
\AA_H([z,w]) = - \int w^*\omega + \int H(t,z(t)) \, dt. \tag 1.12
$$
\endroster
Because our $X_f$ is the negative of $X_f$ in this convention, the
action functional is the one for the Hamiltonian $-H$ in our
convention.  While our convention makes the positive Morse
gradient flow correspond to the negative Cauchy-Riemann flow, the
other convention keeps the same direction. The reason why we keep
our convention is that we would like to keep the definition of the
action functional the same as the classical Hamilton's functional
$$
\int pdq - H\, dt \tag 1.13
$$
on the phase space and to make the {\it negative} gradient flow of
the action functional for the zero Hamiltonian become the
pseudo-holomorphic equation.

It appears that the origin of the two different conventions is the
choice of the convention on how one defines the canonical
symplectic form on the cotangent bundle $T^*N$ or in the classical
phase space: If we set the canonical Liouville form
$$
\theta = \sum_{i} p_idq^i
$$
for the canonical coordinates $q^1, \cdots, q^n, p_1, \cdots, p_n$
of $T^*N$, we take the standard symplectic form to be
$$
\omega_0 = - d\theta = \sum dq^i \wedge dp_i
$$
while the people using the other convention (see e.g., [Po3])
take
$$
\omega_0 = d\theta = \sum dp_i \wedge dq^i.
$$
As a consequence, the action functional (1.12) in the other
convention is the {\it negative} of the classical Hamilton's
functional (1.13). It seems that there is not a single convention
that makes everybody happy and hence one has to live with some
nuisance in this matter one way or the other.

\head{\bf \S 2. The action functional and the action spectrum}
\endhead
Let $(M,\omega)$ be any compact symplectic manifold.
and $\Omega_0(M)$ be the set of contractible loops and
$\widetilde\Omega_0(M)$ be its the covering space mentioned before.
We will always consider {\it normalized} functions $f: M \to \R$
by
$$
\int_M f\, d\mu = 0 \tag 2.1
$$
where $d\mu$ is the Liouville measure of $(M,\omega)$.

When a periodic normalized Hamiltonian $H:M \times
(\R/\Z) \to \R$ is given, we consider the action functional $\AA_H:
\widetilde \Omega(M) \to \R$ by
$$
\AA_H([\gamma,w])= -\int w^*\omega - \int H(\gamma(t),t)dt
$$
We denote by $\text{Per}(H)$ the set of periodic orbits of $X_H$.
\medskip

\definition{Definition 2.1}  We define the {\it
action spectrum} of $H$, denoted as $\hbox{\rm Spec}(H) \subset
\R$, by
$$
\hbox{\rm Spec}(H) := \{\AA_H(z,w)\in \R ~|~ [z,w] \in
\widetilde\Omega_0(M), z\in \text {Per}(H) \},
$$
i.e., the set of critical values of $\AA_H: \widetilde\Omega(M)
\to \R$. For each given $z \in \text {Per}(H)$, we denote
$$
\hbox{\rm Spec}(H;z) = \{\AA_H(z,w)\in \R ~|~ (z,w) \in
\pi^{-1}(z) \}.
$$
\enddefinition

Note that $\text {Spec}(H;z)$ is a principal homogeneous space
modelled by the period group of $(M,\omega)$
$$
\Gamma_\omega := \{ \omega(A)~|~ A \in \pi_2(M)\} = \omega(\Gamma)
\subset \R
$$
and
$$
\hbox{\rm Spec}(H)= \cup_{z \in \text {Per}(H)}\text {Spec}(H;z).
$$

The following was proven in [Oh4].

\proclaim\nofrills{Lemma 2.2.}~ For any closed symplectic manifold
$(M,\omega)$ and for any smooth Hamiltonian $H$, $\hbox{\rm
Spec}(H)$ is a measure zero subset of $\R$ for any $H$.
\endproclaim

For given $\phi \in {\cal H}am(M,\omega)$, we denote $F
\mapsto \phi$ if $\phi^1_F = \phi$, and denote
$$
\HH(\phi) = \{ F ~|~ F \mapsto \phi \}.
$$
We say that two Hamiltonians $F$ and $K$ are equivalent and denote
$F\sim K$ if they
are connected by one parameter family of Hamiltonians
$\{F^s\}_{0\leq s\leq 1}$ such that $F^s \mapsto \phi$ for all $s
\in [0,1]$. We write $[F]$ for the equivalence class of $F$. Then
the universal covering space $\widetilde{{\cal  H}am}(M,\omega)$
of ${\cal  H }am(M,\omega)$ is realized by the set of such
equivalence classes. Note that the group $G:= \Omega({\cal  H
}am(M,\omega),id)$ of based loops
naturally acts on the loop space $\Omega(M)$ by
$$
(h\cdot \gamma) (t) = h(t)(\gamma(t))
$$
where $h \in \Omega({\cal H}am (M,\omega))$ and $\gamma \in
\Omega(M)$. An interesting consequence of Arnold's conjecture is
that this action maps $\Omega_0(M)$ to itself (see e.g., [Lemma
2.2, Se]). Seidel [Lemma 2.4, Se] proves that this action
 can be lifted to $\widetilde\Omega_0(M)$. The set of
lifts $(h,\widetilde h)$ forms a covering group $\widetilde G \to G$
$$
\widetilde G \subset G \times Homeo(\widetilde \Omega_0(M))
$$
whose fiber is isomorphic to $\Gamma$.
Seidel relates the lifting $(h,\widetilde h)$ of
$h: S^1 \to \HH am(M,\omega)$ to a section of the
Hamiltonian bundle associated to the loop $h$ (see \S 2 [Se]).

When a Hamiltonian $H$ generating the loop $h$ is
given, the assignment
$$
z \mapsto h\cdot z
$$
provides a natural one-one correspondence
$$
h: \text{Per}(F) \mapsto \text{Per}(H\# F)
\tag 2.2
$$
where $H\# F = H + F\circ (\phi_H^t)^{-1}$. Let $F, \, K$ be
normalized Hamiltonians with $F, K \mapsto \phi$ and $H$ be the
Hamiltonian such that $K = H \#F$, and $f_t, \, g_t$ and $h_t$ be
the corresponding Hamiltonian paths as above. In particular the
path $h = \{h_t\}_{0\leq t \leq 1}$ defines a loop. We also denote
the corresponding action of $h$ on $\Omega_0(M)$ by $h$. Let
$\widetilde h$ be any lift of $h$ to $\text{Homeo}(\widetilde
\Omega_0(M))$. Then a straightforward calculation shows (see
[Oh5])
$$
\widetilde h^*(d\AA_F) = d\AA_K \tag 2.3
$$
as a one-form on $\widetilde \Omega_0(M)$. In particular since
$\widetilde\Omega_0(M)$ is connected, we have
$$
\widetilde h^*(\AA_F) - \AA_K = C(F,K, \widetilde h)
\tag 2.4
$$
where $C= C(F,K, \widetilde h)$ is a constant a priori depending on
$F, K, \widetilde h$.

\proclaim{Theorem 2.3 [Theorem II, Oh5]}  Let $h$ be the loop as
above and $\widetilde h$ be a lift. Then the constant $C(F,K,
\widetilde h)$ in (2.4) depends only on the homotopy class
$[h,\widetilde h] \in \pi_0(\widetilde G)$.
In particular if $F\sim K$, we have $\AA_F\circ \widetilde h = \AA_K$
and hence
$$
\text{\rm Spec }F = \text{\rm Spec } K $$ as a subset of $\R$. For
any $\widetilde \phi \in \widetilde{\HH am}(M,\omega)$, we define
$$
\text{\rm Spec }(\widetilde \phi) := \text{\rm Spec }F
$$
for a (and so any) normalized Hamiltonian $F$ with
$[\phi,F] = \widetilde \phi$.
\endproclaim

\definition{Definition 2.4 [Action Spectrum Bundle]}
We define the action spectrum bundle of $(M,\omega)$ by
$$
\frak{Spec}(M,\omega) =  \bigcup_{\widetilde \phi \in
\widetilde \HH am(M,\omega)}\frak{Spec}_{\widetilde \phi}(M,\omega)
\subset \widetilde{\HH am}(M,\omega) \times \R
$$
where
$$
\frak{Spec}_{\widetilde \phi}(M,\omega)
= \{\AA_F([z,w]) \mid d\AA_F([z,w]) = 0 , \quad \widetilde \phi
= [F ] \, \}  \subset \R
$$
and denote by $\pi: \frak{Spec}(M,\omega)
\to \widetilde{\HH am}(M,\omega)$ the natural projection.
\enddefinition

\head{\bf \S 3. Quantum cohomology in the chain level}
\endhead

We first recall the definition of the quantum cohomology ring
$QH^*(M)$. As a module, it is defined as
$$
QH^*(M) = H^*(M,\Q) \otimes \Lambda^\uparrow_\omega
$$
where $\Lambda_\omega^\uparrow$ is the
(upward) Novikov ring
$$
\Lambda_\omega^\uparrow = \Big\{
\sum_{A \in \Gamma} a_A q^{-A} \mid
a_A \in \Q, \, \# \{ A \mid a_i\neq 0, \, \int_{-A} \omega < \lambda \}
< \infty, \, \forall \lambda \in \R\Big \}.
$$
Due to the finiteness assumption on the Novikov ring, we have the
natural (upward) valuation $v: QH^*(M) \to \R$ defined by
$$
v(\sum_{A \in \Gamma_\omega} a_A q^{-A}) = \min\{\omega(-A): a_A
\neq 0\} \tag 3.1
$$
which satisfies that for any $a, \ b \in QH^*(M)$
$$
v(a+b) \geq \min\{v(a), v(b)\}.
$$
\definition{Definition 3.1}
For each homogeneous element
$$
a = \Sigma_{A \in \Gamma} a_A q^{-A} \in QH^k(M),
\quad a_A \in H^*(M,\Q)
\tag 3.2
$$
of degree $k$,  we also call $v(a)$ the {\it level} of $a$ and the
corresponding term in the sum the {\it leading order term} of $a$
and denote by $\sigma(a)$. Note that the leading order term
$\sigma(a)$ of a {\it homogeneous} element $a$ is unique among the
summands in the sum by the definition (1.4) of $\Gamma$.
\enddefinition

The product on $QH^*(M)$ is defined by the usual quantum cup product, which we
denote by ``$\cdot$'' and which preserves the grading, i.e, satisfies
$$
QH^k(M) \times QH^\ell(M) \to QH^{k+\ell}(M).
$$
Often the homological version of the quantum cohomology is also useful,
sometimes called the quantum homology, which is defined by
$$
QH_*(M) = H_*(M) \otimes \Lambda_\omega^\downarrow
$$
where $\Lambda^\downarrow_\omega$ is the (downward) Novikov ring
$$
\Lambda_\omega^\downarrow = \{
\sum_{B_j \in \Gamma} b_j q^{B_j} \mid
b_j \in \Q, \, \# \{ B_j \mid b_j\neq 0, \, \int_{B_j} \omega > \lambda \}
< \infty, \forall \lambda \in \R\}.
$$

We define the corresponding (downward) valuation by
$$
v(\sum_{B \in \Gamma} a_B q^{B}) = \max\{\omega(B): a_B \neq 0\}
\tag 3.3
$$
which satisfies that for $f,\, g \in QH_*(M)$
$$
v(f + g) \leq \max\{v(f), v(g) \}.
$$

We like to point out that the summand in
$\Lambda_\omega^\downarrow$ is written as $b_B q^B$ while the one
in $\Lambda_\omega^\uparrow$ as $a_A q^{-A}$ with the minus sign.
This is because we want to clearly show which one we use.
Obviously $-v$ in (3.1) and $v$ in (3.3) satisfy the axiom of
non-Archimedean norm which induce a topology on $QH^*(M)$ and
$QH_*(M)$ respectively. In each case the finiteness assumption in
the definition of the Novikov ring allows us to numerate the
non-zero summands in each given Novikov chain (3.2) so that
$$
\lambda_1 > \lambda_2 > \cdots > \lambda_j > \cdots
$$
with $\lambda_j = \omega(B_j)$ or $\omega(A_j)$.

Since the downward Novikov ring appears mostly in this paper, we will just use
$\Lambda_\omega$ or $\Lambda$ for $\Lambda^\downarrow_\omega$, unless absolutely
necessary to emphasize the direction of the Novikov ring.
We define the {level} and the {leading order term} of $b \in QH_*(M)$
similarly as in Definition 3.1 by changing the role of upward and downward
Novikov rings. We have a canonical isomorphism
$$
\flat: QH^*(M) \to QH_*(M); \quad \sum a_i q^{-A_i} \to
\sum PD(a_i) q^{A_i}
$$
and its inverse
$$
\sharp: QH_*(M) \to QH^*(M); \quad \sum b_j q^{B_j} \to \sum
PD(b_j) q^{-B_j}.
$$
We denote by $a^\flat$ and $b^\#$ the images under these maps.

There exists the canonical non-degenerate pairing
$$
\langle \cdot, \cdot \rangle : QH^*(M) \otimes QH_*(M) \to \Q
$$
defined by
$$
\langle \sum a_i q^{-A_i}, \sum b_j q^{B_j} \rangle =
\sum (a_i,b_j) \delta_{A_iB_j}
\tag 3.4
$$
where $\delta_{A_iB_j}$ is the delta-function and $(a_i,b_j)$ is the
canonical pairing between $H^*(M,\Q)$ and $H_*(M,\Q)$.
Note that this sum is always finite by the finiteness condition
in the definitions of $QH^*(M)$ and $QH_*(M)$ and so is well-defined.
This is equivalent to the Frobenius pairing in the quantum cohomology
ring. However we would like to emphasize that the dual vector space
$(QH_*(M))^*$ of $QH_*(M)$ is {\it not} isomorphic to $QH^*(M)$ even
as a $\Q$-vector space. Rather the above pairing induces an injection
$$
QH^*(M) \hookrightarrow (QH_*(M))^*
$$
whose images lie in the set of {\it continuous} linear functionals
on $QH_*(M)$ with respect to the topology induced by the valuation
$v$. (3.3) on $QH_*(M)$. We refer to [Br] for a good introduction
to non-Archimedean analytic geometry. In fact, the description of
the standard quantum cohomology in the literature is {\it not}
really a `cohomology' but a `homology' in that it never uses
linear functionals in its definition. To keep our exposition
consistent with the standard literature in the Gromov-Witten
invariants and the quantum cohomology, we prefer to call them the
quantum cohomology rather than the quantum homology as some
authors did (e.g., [Se]) in the symplectic geometry community. In
Appendix, we will introduce a {\it genuinely cohomological}
version of quantum cohomology which we call {\it continuous
quantum cohomology} using the continuous linear functionals on the
{\it quantum chain complex} below with respect to the topology
induced by the valuation $v$.

Let $(C_*, \partial)$ be any chain complex on $M$ whose homology
is the singular homology $H_*(M)$. One may take for $C_*$ the
usual singular chain complex or the Morse chain complex of a Morse
function $f: M \to \R$, $(C_*(-\e f), \partial_{-\e f})$ for some
sufficiently small $\e > 0$. However since we need to take a
non-degenerate pairing in the chain level, we should use a model
which is {\it finitely generated}. We will always prefer to use
the Morse homology complex because it is finitely generated and
avoids some technical issue related to singular degeneration
problem of the type studied in [FOh1,2]. The negative sign in
$(C_*(-\e f),
\partial_{-\e f})$ is put to make the correspondence between the
Morse homology and the Floer homology consistent with our
conventions of the Hamiltonian vector field (1.2) and the action
functional (1.6). In our conventions, solutions of negative
gradient of $-\e f$ correspond to ones for the negative gradient
flow of the action functional $\AA_{\e f}$. We denote by
$$
(C^*(-\e f), \delta_{-\e f})
$$
the corresponding cochain complex, i.e,
$$
C^k := \text{Hom}(C_k, \Q), \quad \delta_{-\e f}
= \partial_{-\e f}^*.
$$

Now we extend the complex $(C_*(-\e f), \partial_{-\e f})$
to the {\it quantum chain complex}, denoted by
$$
(CQ_*(-\e f), \partial_Q)
$$
$$
CQ_*(-\e f) : = C_*(-\e f) \otimes \Lambda_\omega, \quad
\partial_Q: = \partial_{-\e f} \otimes \Lambda_\omega.
 \tag 3.5
$$
This coincides with the Floer complex $(CF_*(\e f), \part)$ as
a chain complex if $\e$ is sufficiently small.
Similarly we define the quantum cochain complex
$(CQ^*(-\e f),\delta^Q)$ by changing the downward Novikov ring to
the upward one. In other words, we define
$$
CQ^*(-\e f): = CM_{2n -*}(-\e f)\otimes \Lambda^\uparrow,
\quad \delta^Q: = \part_{\e f}\otimes
\Lambda^\uparrow_\omega.
$$
Again we would like to emphasize that $CQ^*(-\e f)$ is {\it not}
isomorphic to the dual space of $CQ_*(-\e f)$ as a $\Q$-vector space.
We refer to Appendix for some further discussion on this issue.

It is well-known that the corresponding homology of
this complex is independent of the choice $f$ and
isomorphic to the above quantum cohomology (resp. the
quantum homology) as a ring (see [PSS],\, [LT2], \, [Lu] for its proof).
This isomorphism however plays no significant role in the current paper,
except for the purpose of bookeeping  the family of invariants $\rho(H;a)$
that we associate to each quantum cohomology class $a \in QH^*(M)$
later (See section 5.1 for more explanation on this point).
To emphasize the role of the Morse function in the level of
complex, we denote the corresponding homology by
$HQ^*(-\e f) \cong QH^*(M)$.

With these definitions, we have the obvious non-degenerate pairing
$$
CQ^*(-\e f) \otimes CQ_*(-\e f) \to \Q \tag 3.6
$$
in the chain level which induces the pairing (3.4) above in homology.

We now choose a generic Morse function $f$. Then
for any given homotopy
$\HH=\{H^s\}_{s \in [0,1]}$ with $H^0 = \e f$ and $H^1 = H$,
we denote by
$$
h_\HH: CQ_*(-\e f) = CF_{*-n}(\e f) \to CF_{*-n}(H) \tag 3.7
$$
the standard Floer chain map from $\e f$ to $H$ via the
homotopy $\HH$.  This induces a homomorphism
$$
h_\HH: HQ_*(-\e f) \cong HF_{*-n}(\e f) \to HF_{*-n}(H). \tag 3.8
$$
Although (3.7) depends on the choice $\HH$, (3.8) is canonical, i.e,
does not depend on the homotopy $\HH$. One confusing point in this
isomorphism is the issue of grading. See the next section for
a review of the construction of this chain map and the issue of
grading of $HF_*(H)$.

\head{\bf \S 4. Filtered Floer homology}
\endhead

For each given generic non-degenerate $H:S^1 \times M \to \R $,
we consider the free $\Q$ vector space over
$$
\text{Crit}\AA_H = \{[z,w]\in \widetilde\Omega_0(M) ~|~ z \in
\text{Per}(H)\}. \tag 4.1
$$
To be able to define the Floer boundary operator correctly, we
need to complete this vector space downward with respect to the
real filtration provided by the action $\AA_H([z,w])$ of the
element $[z, w]$ of (4.1). More precisely,
\medskip

\definition {Definition 4.1} We call the formal sum
$$
\beta = \sum _{[z, w] \in \text{Crit}\AA_H} a_{[z, w]} [z,w], \,
a_{[z,w]} \in \Q \tag 4.2
$$
a {\it Novikov chain} if there are
only finitely many non-zero terms in the expression (4.2) above
any given level of the action. We denote by $CF_*(H)$
the set of Novikov chains. We call those $[z,w]$ with $a_{[z,w]} \neq 0$
{\it generators} of the chain $\beta$ and just denote as
$$
[z,w] \in \beta
$$
in that case. Note that $CF_*(H)$ is a graded
$\Q$-vector space which is infinite dimensional in general,
unless $\pi_2(M) = 0$.
\enddefinition

As in [Oh4], we introduce the following notion which is a crucial concept
for the mini-max argument we carry out later.

\definition{Definition 4.2}~  Let $\beta$ be a Novikov chain in
$CF_*(H)$. We define the {\it level} of the cycle
$\beta$ and denote by
$$
\lambda_H(\beta) =\max_{[z,w]} \{\AA_H([z,w]) ~|~a_{[z,w]}  \neq
0\, \, \text{in }\, (4.2) \}
$$
if $\beta \neq 0$, and just put $\lambda_H(0) = - \infty$ as
usual. We call the unique critical point $[z,w]$ that realizes the
maximum value $\lambda_H(\beta)$ the {\it peak} of the cycle
$\beta$, and denote it by $pk(\beta)$.
\enddefinition

We briefly review construction of basic operators in the Floer
homology theory [Fl]. Let $J = \{J_t\}_{0\leq t \leq 1}$ be a
periodic family of compatible almost complex structure on
$(M,\omega)$.

For each given pair $(J, H)$, we define the boundary operator
$$
\part:CF_*(H) \to CF_*(H)
$$
considering the perturbed Cauchy-Riemann equation
$$
\cases
\frac{\part u}{\part \tau} + J\Big(\frac{\part u}{\part t}
- X_H(u)\Big) = 0\\
\lim_{\tau \to -\infty}u(\tau) = z^-,  \lim_{\tau \to
\infty}u(\tau) = z^+ \\
\endcases
\tag 4.3
$$
This equation, when lifted to $\widetilde \Omega_0(M)$, defines
nothing but the {\it negative} gradient flow of $\AA_H$ with
respect to the $L^2$-metric on $\widetilde \Omega_0(M)$ induced by
the family of metrics on $M$
$$
g_{J_t} = (\cdot, \cdot)_{J_t}: = \omega(\cdot, J_t\cdot):
$$
This $L^2$-metric is defined by
$$
\langle \xi, \eta \rangle_J: = \int_0^1 \ (\xi, \eta)_{J_t}\, dt.
$$
We will also denote
$$
\|v\|^2_{J_0} = (v,v)_{J_0} = \omega(v, J_0 v) \tag 4.4
$$
for $v \in TM$.

For each given
$[z^-,w^-]$ and $[z^+,w^+]$, we define the moduli space
$$
\MM(H,J; [z^-,w^-],[z^+,w^+])
$$
of solutions $u$ of (4.3) wit finite energy
$$
E_J(u) = \frac{1}{2} \int_{\R \times S^1}
\Big(\Big|\dudtau\Big|_{J_t}^2 + \Big|\dudt - X_H(u)\Big|_{J_t}^2
\Big)dt\,d\tau < \infty
$$
and  satisfying
$$
w^-\# u \sim w^+. \tag 4.5
$$
$\part$ has degree $-1$ and satisfies $\part\circ \part = 0$.

When we are given a family $(j, \HH)$ with $\HH = \{H^s\}_{0\leq s
\leq 1}$ and $j = \{J^s\}_{0\leq s \leq 1}$, the chain
homomorphism
$$
h_{(j,\HH)}: CF_*(H^0) \to CF_*(H^1)
$$
is defined by the non-autonomous equation
$$
\cases \frac{\part u}{\part \tau} +
J^{\rho_1(\tau)}\Big(\frac{\part u}{\part t}
- X_{H^{\rho_2(\tau)}}(u)\Big) = 0\\
\lim_{\tau \to -\infty}u(\tau) = z^-,  \lim_{\tau \to
\infty}u(\tau) = z^+
\endcases
\tag 4.6
$$
also with the condition (4.5).
Here $\rho_i, \, i= 1,2$ is the cut-off functions of the type $\rho :\R \to
[0,1]$,
$$
\align
\rho(\tau) & = \cases 0 \, \quad \text {for $\tau \leq -R$}\\
                    1 \, \quad \text {for $\tau \geq R$}
                    \endcases \\
\rho^\prime(\tau) & \geq 0
\endalign
$$
for some $R > 0$.  $h_{(j,\HH)}$ has degree 0 and satisfies
$$
\part_{(J^1,H^1)} \circ h_{(j,\HH)} = h_{(j,\HH)} \circ
\part_{(J^0,H^0)}.
$$
Two such chain maps for different homotopies $(j^1,\HH^1)$
and $(j^2, \HH^2)$ connecting the same end points
are also known to be chain homotopic [Fl2].

Finally when we are given a homotopy $(\overline j, \overline
\HH)$ of homotopies with $\overline j = \{j_\kappa\}$,
$\overline\HH = \{\HH_\kappa\}$, consideration of the
parameterized version of (4.6) for $ 0 \leq \kappa \leq 1$ defines
the chain homotopy map
$$
H_{\overline\HH} :CF_*(H^0) \to CF_*(H^1)
\tag 4.7
$$
which has degree $+1$ and satisfies
$$
h_{(j_1, \HH_1)} - h_{(j_0,\HH_0)} = \part_{(J^1,H^1)} \circ
H_{\overline\HH} + H_{\overline\HH} \circ \part_{(J^0,H^0)}.
\tag 4.8
$$
By now, construction of these maps using these moduli spaces has
been completed with rational coefficients (See [FOn], [LT1] and
[Ru]) using the techniques of virtual moduli cycles.
We will suppress this advanced machinery from our presentation
throughout the paper. The main stream of the proof is
independent of this machinery except that it is implicitly
needed to prove that various moduli spaces we use are non-empty.
Therefore we do
not explicitly mention these technicalities in the main body of
the paper until \S 8, unless it is absolutely necessary.
In \S 8, we will provide justification of this in the general case
all at once.

The following upper estimate of the action change can be proven by
the same argument as that of the proof of [Ch], [Oh1,4]. We would
like to emphasize that in general there does {\it not} exist a
lower estimate of this type. The upper estimate is just one
manifestation of the `positivity' phenomenon in symplectic
topology through the existence of pseudo-holomorphic curves that
was first discovered by Gromov [Gr].  On the other hand, the existence of
lower estimate is closely tied to some nontrivial homological
property of (Floer) cycles, and best formulated in terms of Floer
cycles instead of individual critical points $[z,w]$ {\it for the
nondegenerate Hamiltonians}.  However, we would like to point out
that the equations (4.3), (4.6) themselves or the numerical estimate of
the action changes for solutions $u$ with finite energy can be studied for any $H$
or $(\HH,j)$ which are not necessarily non-degenerate or generic,
although the Floer complex or the operators may not be defined for
such choices.

\proclaim{Proposition 4.3} Let $H, K$ be any Hamiltonian not
necessarily non-degenerate and $j = \{J^s\}_{s \in [0,1]}$ be any
given homotopy and $\HH^{lin} = \{H^s\}_{0\leq s\leq 1}$ be the
linear homotopy $H^s = (1-s)H + sF$. Suppose that (4.6) has a
solution satisfying (4.5). Then we have the identity
$$
\align \AA_F([z^+,w^+]) & - \AA_H([z^-,w^-]) \\
& = - \int \Big|\dudtau \Big|_{J^{\rho_1(\tau)}}^2 -
\int_{-\infty}^\infty \rho'(\tau)\Big(F(t,u(\tau,t)) -
H(t,u(\tau,t))\Big)
\, dt\,d\tau  \tag 4.9\\
& \leq - \int \Big|\dudtau \Big|_{J^{\rho_1(\tau)}}^2 + \int_0^1
-\min_{x \in M} (F_t - H_t) \, dt \tag 4.10\\
& \leq  \int_0^1 -\min_{x \in M} (F_t - H_t) \, dt \tag 4.11
\endalign
$$
\endproclaim

By considering the case $K=H$, we immediately have

\proclaim{Corollary 4.4 }~ For a
fixed $H$ and for a given one parameter family $j =
\{J^s\}_{s \in [0,1]}$, let $u$ be as in Proposition 4.3. Then
we have
$$
\AA_H([z^+,w^+]) - \AA_H([z^-,w^-]) = - \int \Big| \dudtau
\Big|_{J^{\rho_1(\tau)}}^2 \leq 0. \tag 4.12
$$
\endproclaim

\definition{Remark 4.5}
We would like to remark that similar calculation proves that there
is also an uniform upper bound $C(j,\HH)$ for the chain map over
general homotopy $(j,\HH)$ or for the chain homotopy maps (4.7).
In this case, the identity (4.9) is replaced by
$$
\align
\AA_F([z^+,w^+]) & - \AA_H([z^-,w^-]) \\
& = - \int \Big|\dudtau \Big|_{J^{\rho_1(\tau)}}^2 -
\int_{-\infty}^\infty \rho'(\tau)\Big( \frac{\part H^s}{\part s}\Big|_{s = \rho(\tau)}
(t, u(\tau,t))\Big)
\, dt\,d\tau \\
& \leq - \int \Big|\dudtau \Big|_{J^{\rho_1(\tau)}}^2 + \int_0^1
-\min_{x \in M} \Big( \frac{\part H^s}{\part s}\Big|_{s =
\rho(\tau)}\Big)\, dt \\
& \leq \int_0^1 -\min_{x \in M} \Big(
\frac{\part H_t^s}{\part s} \Big)\, dt \tag 4.13
\endalign
$$
This upper estimate is also crucial for the construction of these maps.
This upper estimate
depends on the choice of homotopy $(j,\HH)$ and is related to the
curvature estimates of the relevant Hamiltonian fibration (see
[Po2], [En1]).
\enddefinition

Now we recall that $CF_*(H)$ is also a
$\Lambda$-module: each $A \in \Gamma$ acts on
${\text Crit}\AA_H$ and so on $CF_*(H)$ by ``gluing a
sphere''
$$
A: [z,w] \to [z, w\# A].
$$
Then $\partial$ is $\Lambda$-linear and induces the
standard Floer homology $HF_*(H;\Lambda)$ with
$\Lambda$ as its coefficients (see [HoS] for a detailed
discussion on the Novikov ring and  on the Floer complex as
a $\Lambda$-module).  However the action
does {\it not} preserve the filtration we defined above.
Whenever we talk about filtration, we will always presume that
the relevant coefficient ring is $\Q$.

For a given nondegenerate $H$ and an $\lambda \in \R \setminus
\hbox{Spec}(H)$, we define the relative chain group
$$
CF_k^\lambda(H): = \{ \beta \in CF_k(H) \mid \lambda_H(\beta) <
\lambda\}.
$$
Corollary 4.4 impies that between the two chain complexes
$(CF_k(H), \part_{(H,J)})$ and   $(CF_k(H), \part_{(H,J')}$,
there is a canonical filtration preserving chain isomorphism
$$
h_{(j, H)}: (CF_k(H), \part_{(H,J)}) \to (CF_k(H), \part_{(H,J')})
$$
where $j$ is any homotopy from $J$ and $J'$, and $\HH\equiv H$ is
the constant homotopy of $H$.  Therefore from now on, we suppress
$J$-dependence on the Floer homology in our exposition unless it
is absolutely necessary.

For each given pair of real numbers $\lambda,  \mu \in \R \setminus
\hbox{Spec}(H)$ with $\lambda < \mu$, we define
$$
CF_*^{(\lambda,\mu]}(H): = CF^\mu(H) / CF^\lambda(H).
$$
Then for each triple $\lambda < \mu < \nu $ where $\lambda=-\infty$
or $\nu=\infty$ are allowed,
we have the short-exact sequence of the complex of graded $\Q$ vector
spaces
$$
0 \to CF^{(\lambda, \mu]}_k(H) \to CF^{(\lambda, \nu]}_k(H) \to
CF^{(\mu, \nu]}_k(H) \to 0
$$
for each $k \in \Z$. This then induces the long exact sequence
of graded modules
$$
\cdots \to HF^{(\lambda, \mu]}_k(H) \to HF^{(\lambda, \nu]}_k(H) \to
HF^{(\mu, \nu]}_k(H) \to HF^{(\lambda, \mu]}_{k-1}(H) \to
\cdots
$$
whenever the relevant Floer homology groups are defined.

We close this section by fixing our grading convention for
$HF_*(H)$. This convention is the analog to the one we use in
[Oh1,2] in the context of Lagrangian submanifolds. We first recall
that solutions of the {\it negative} gradient flow equation of
$-f$, (i.e., of the {\it positive} gradient flow of $f$
$$
\dot\chi - \text{grad } f(\chi) = 0
$$
corresponds to the {\it negative} gradient flow of the action functional
$\AA_{\e f}$). This gives rise to
the relation between the Morse indices $\mu_{-\e f}^{Morse}(p)$
and the Conley-Zehnder indices $\mu_{CZ}([p,\widehat p];\e f)$
(see [Lemma 7.2, SZ] but with some care about the different convention of
the Hamiltonian vector field. Their definition of $X_H$ is $-X_H$
in our convention):
$$
\mu_{CZ}([p,\widehat p];\e f)  = \mu_{-\e f}^{Morse}(p) -n \tag
4.14
$$
in our convention. On the other hand, obviously we have
$$
\mu_{-\e f}^{Morse}(p) - n = (2n-\mu_{\e f}^{Morse}(p)) -n = n -
\mu_{\e f}^{Morse}(p)
$$
We will always grade $HF_*(H)$ by the Conley Zehnder index
$$
k = \mu_H([z,w]):= \mu_{CZ}([z,w];H). \tag 4.15
$$
This grading convention makes the degree $k$ of $[q,\widehat q]$ in
$CF_k(\e f)$ coincides with the Morse index of $q$ of $\e f$ for each
$q \in \text{Crit}\e f$. Recalling that we chose the Morse complex
$$
CM_*(-\e f) \otimes \Lambda^\downarrow
$$
for the quantum chain complex $CQ_*(-\e f)$, it also coincides with the standard
grading of the quantum cohomology via the map
$$
\flat: QH^k(M) \to QH_{2n-k}(M).
$$
Form now on, we will just denote by $\mu_H([z,w])$ the Conley-Zehnder
index of $[z,w]$ for the Hamiltonian $H$.
Under this grading, we have the following grading preserving isomorphism
$$
QH^{n-k}(M) \to QH_{n+k}(M) \cong HQ_{n+k}(-\e f) \to HF_k(\e f) \to
HF_k(H). \tag 4.16
$$
We will also show in \S 6 that this grading convention makes the
pants product, denoted by $*$, has the degree $-n$:
$$
*: HF_k(H) \otimes HF_\ell(K) \to HF_{(k+\ell)-n}(H\# K)
\tag 4.17
$$
which will be compatible with the degree preserving quantum product
$$
\cdot : QH^a(M) \otimes QH^b(M) \to QH^{a+b}(M).
$$

\head{\bf \S 5. Construction of the spectral invariants of
Hamiltonian functions}
\endhead

In this section, we associate some homologically essential
critical values of the action functional $\AA_H$ to each
Hamiltonian functions $H$ and quantum cohomology class $a$, and
call them the {\it spectral invariants} of $H$. We denote this
assignment by
$$
\rho: C^\infty_m([0,1]\times M) \times QH^*(M) \to \R
$$
as described in the introduction of this paper. Before launching our
construction, some overview of our construction of spectral invariants
is necessary.

\medskip
\n{\it 5.1. Overview of the construction}
\smallskip

We recall the canonical isomorphism
$$
h_{\alpha\beta}: HF_*(H_\alpha) \to HF_*(H_\beta)
$$
which satisfies the composition law
$$
h_{\alpha\gamma} = h_{\alpha\beta}\circ h_{\beta\gamma}.
$$
We denote by $HF_*(M)$ the corresponding model $\Q$-vector space.
We also note that $HF_*(H)$ is induced by the filtered chain complex
$(CF_*^\lambda (H), \part)$ where
$$
CF_*^\lambda(H) = \operatorname{span}_\Q
\{\alpha \in CF_*(H) \mid \lambda_H(\alpha)  \leq \lambda \}
$$
i.e., the sub-complex generated by the critical points $[z,w] \in
\text{Crit}\AA_H$ with
$$
\AA_H([z,w]) \leq \lambda.
$$
Then there exists a canonical inclusion
$$
i_\lambda : CF_*^{\lambda}(H) \to CF_*^\infty(H):=CF_*(H)
$$
which induces a natural homomorphism $i_\lambda: HF_*^\lambda(H) \to
HF_*(H)$.
For each given element $\ell \in FH_*(M)$ and Hamiltonian $H$, we represent
the class $\ell$ by a Novikov cycle $\alpha$ of $H$ and measure its level
$\lambda_H(\alpha)$ and define
$$
\rho(H;\ell):= \inf\{ \lambda \in \R\mid \ell \in \text{Im }i_\lambda\}
$$
or equivalently
$$
\rho(H;\ell): = \inf_{\alpha; i_\lambda[\alpha] = \ell} \lambda_H(\alpha).
$$
The crucial task then is to prove that for each (homogeneous)
element $\ell \neq 0$, the value $\rho(H;\ell)$ is finite, i.e,
``the cycle $\alpha$ is linked and cannot be pushed away to
infinity by the negative gradient flow of the action functional''.
In the classical critical point theory (see [BnR] for example),
this property of semi-infinite cycles is called the {\it linking }
property. We like to point out that there is no manifest way to
see the linking property or the criticality of the mini-max value
$\rho(H;\ell)$ out of the definition itself.

We will prove this finiteness first for the Hamiltonian $\e f$
where $f$ is a Morse function and $\e$ is sufficiently small. This
finiteness strongly relies on the facts that the Floer boundary
operator $\part_{\e f}$ in this case has the form
$$
\part_{\e f} = \part_{-\e f}^{Morse} \otimes \Lambda_\omega
\tag 5.1
$$
i.e, ``there is no quantum contribution on the Floer boundary
operator'', and that the classical Morse theory proves that
$\part_{-\e f}^{Morse}$ cannot push down the level of a
non-trivial cycle more than $-\e \max f$ (see [Oh4]).

Once we prove the finiteness for $\e f$, then we
consider the general nondegenerate Hamiltonian $H$.
We compare the cycles in $CF_*(H)$ and the transferred cycles
in $CF_*(\e f)$ by the chain
map $h_{\HH}^{-1}: CF_*(H) \to CF_*(\e f)$ where $\HH$ is a homotopy
connecting $\e f$ and $H$. The change of the level then can be
measured by judicious use of (4.7) and Remark 4.5 which will prove
the finiteness for any $H$.

After we prove finiteness of $\rho(H;a)$ for general $H$,
we study the continuity property of $\rho(H;a)$
under the change of $H$.
This will be done, via the equation (4.6),  considering the level
change between arbitrary pair $(H,K)$.

Finally we  prove that the limit
$$
\lim_{\e \to 0} \rho(\e f;\ell)
$$
exists and is independent of the choice of Morse function $f$.
If the Floer homology class $\ell$ is identified with $a^\flat$
for a quantum cohomology class $a \in QH^*(M)$ under the
PSS-isomorphism [PSS], then this limit is nothing but  the
valuation $v(a)$.

In this procedure, we can avoid considering the `singular limit'
of the `chains' (See the [section 8, Oh8] for some illustration of
the difficulty in studying such limits). We only need to consider
the limit of the values $\rho(H;\ell)$ as $H \to 0$ which is a
much simpler task than considering the limit of chains which
involves highly non-trivial analytical work (we refer to the
forthcoming work [FOh2] for the consideration of this limit in the
chain level).

\medskip
\n{\it 5.2. Finiteness; the linking property of semi-infinite
cycles}
\smallskip

With this overview, we now start with our construction.
We first recall the natural pairing
$$
\langle \cdot, \cdot \rangle: CQ^*(-\e f) \otimes CQ_*(-\e f) \to \Q:
$$
where we have
$$
\align
CQ_k(-\e f) & := (CM_k(-\e f), \part_{-\e f}) \otimes \Lambda^\downarrow\\
CQ^k(-\e f) & := (CM_{2n-k}(\e f), \part_{\e f}) \otimes \Lambda^\uparrow. \\
\endalign
$$
\definition{Remark 5.1} We would like to emphasize that in our definition
$CQ^k(-\e f)$ is not isomorphic to $\text{Hom}_\Q(CQ_k(-\e f), \Q)$ in general.
However there is a natural homomorphism
$$
CQ^k(-\e f) \to \text{Hom}_\Q(CQ_k(-\e f), \Q); \quad a \mapsto
\langle a, \cdot \rangle \tag 5.2
$$
whose image lies in the subset of continuous linear functionals
$$
Hom_{cont}(CQ_k(-\e f), \Q): = CQ^k_{cont}(-\e f) \subset Hom_{\Q}
(CQ_k(-\e f), \Q).
$$
See Appendix for more discussions on this aspect. We would like to
emphasize that (5.2) is well-defined because of the choice of
directions of the Novikov rings $\Lambda^\uparrow$ and
$\Lambda^\downarrow$. In general, the map (5.2) is injective but
not an isomorphism. Polterovich [Po4], [EnP] observed that this
point is closely related to certain failure of ``Poincar\'e
duality'' of the Floer homology with Novikov rings as its
coefficients.
\enddefinition

Now we are ready to give the definition of our spectral
invariants. Previously in [Oh4], the author outlined this
construction for the classical cohomology class in $H^*(M) \subset
QH^*(M)$.

\definition{Definition 5.2} Let $H$ be a generic non-degenerate Hamiltonian.
For each given $a \in QH^k(M)\cong HQ^k(-\e f)$, we define
$$
\rho(H,a) = \inf_\alpha \{ \lambda_H(\alpha) \mid [\alpha] =
a^\flat,\, \alpha \in CF_k(H) \}. \tag 5.3
$$
\enddefinition

\proclaim{Theorem 5.3}  Let $0 \neq a \in QH^*(M)$.
\roster \item
Let $H$ be a generic non-degenerate Hamiltonian. Then $\rho(H,a)$ is finite.
\item
For any pair of generic nondegenerate Hamiltonians $H, \, K$, we have the
inequality
$$
\int_0^1 - \max(K-H)\, dt  \leq \rho(K,a) -
\rho(H,a) \leq  \int_0^1 - \min(K-H)\, dt.
\tag 5.4
$$
In particular, the function $H \mapsto \rho(H;a)$ continuously
extends to $C^0_m([0,1]\times M)$.
\endroster
\endproclaim
\demo{Proof} We will prove the finiteness in two steps: first we prove
the finiteness for $\e f$ for sufficiently small $\e >0$ for any
given Morse function $f$, and then prove it for general $H$ using
this finiteness for $\e f$. After then we will prove the
inequality (5.4).

\medskip
\n{\it Step 1: The finiteness of for $\e f$.}  Let $f$ be any
fixed Morse function and fix $\e > 0$ so small that there is no
quantum contribution for the Floer boundary operator $\part_{(\e
f, J_0)}$ for a time independent family $J_t \equiv J_0$ for any
compatible almost complex structure $J_0$, i.e, we have
$$
\part_{(\e f, J_0)} \simeq \part_{-\e f}^{Morse}\otimes
\Lambda_\omega^\downarrow. \tag 5.5
$$
It is well-known ([Fl], [FOn], [LT1]) that this is possible.
Fixing such $\e$ and $J_0$, we just denote
$$
\part_{\e f} = \part_{(\e f,J_0)}.
$$
Then by considering the Morse homology of $- \e f$ with respect to
the Riemannian metric $g_{J_0}= \omega(\cdot, J_0 \cdot)$, we have
the identity
$$
\align QH^*(M) & \cong \ker \part_{\e f}^{Morse}\otimes
\Lambda^\uparrow /\text{Im } \part_{\e f}^{Morse}\otimes
\Lambda^\uparrow =
HM_*(\e f) \otimes \Lambda^\uparrow\\
QH_*(M) & \cong \ker \part_{-\e f}^{Morse}\otimes
\Lambda^\downarrow /\text{Im } \part_{-\e f}^{Morse}\otimes
\Lambda^\downarrow = HM_*(-\e f) \otimes \Lambda^\downarrow.
\endalign
$$
Recalling
$$
CF_k(\e f) \cong CQ_{n+k}(-\e f),
$$
from (5.5), we represent $a^\flat \in QH_{n+k}(M)$ by a Novikov
cycle of $\e f$ where
$$
\alpha = \sum_A a_{p \otimes q^{A}} \, p \otimes q^{A}
$$
with $a_p \in \Q$ and $p \in \text{Crit}_*(-\e f)$ and
$$
n+k = \mu_{\e f}(p\otimes q^A) \tag 5.6
$$
where $\mu_{\e f}(p \otimes q^A)$ is the Conley-Zehnder index of
the element $p\otimes q^{A} = [p, \widehat p \# A]$. We  recall
the general index formula
$$
\mu_H([z,w\otimes A]) = \mu_H([z,w]) + 2c_1(A)
$$
in our convention (see [Oh9] for the proof of this index formula).
Applying this to $H = \e f$, we have obtained
$$
\mu_{\e f}([p, \widehat p\# A]) = \mu_{\e f}([p, \widehat p]) +
2c_1(A).
$$
Combining this with
$$
\mu^{Morse}_{-\e f}(p) = \mu_{\e f}([p, \widehat p]) +n
$$
we derive that (5.6) is equivalent to
$$
\mu^{Morse}_{-\e f}(p) = n+ k - 2c_1(A).
$$
Next we see that $\alpha $ has the level
$$
\lambda_{\e f}(\alpha) = \max\{-\e f(p) - \omega(A) \mid
a_{p \otimes q^{A}} \neq 0 \} \tag 5.7
$$
because $\AA_{\e f}([p, \widehat p \# A]) = -\e f(p) - \omega(A)$.
Now the most crucial point in our construction is to prove
the finiteness
$$
\inf_{[\alpha] = a^\flat}\lambda_{\e f}(\alpha) >
-\infty. \tag 5.8
$$
The following lemma proves this linking property. We first like to
point out that the quantum cohomology class
$$
a = \sum_A a_A q^{-A}
$$
uniquely determines the set
$$
\Gamma(a): =\{ A \in \Gamma \mid a_A \neq 0\}.
$$
By the finiteness condition in the formal power series, we can
enumerate $\Gamma(a)$ so that
$$
\lambda_1 > \lambda_2 > \lambda_3 > \cdots \tag 5.9
$$
without loss of generality. In particular, we have
$$
v(a) = -\omega(A_1)= \lambda_1. \tag 5.10
$$
\proclaim{Lemma 5.4} Let $a \neq QH^k(M)$ and $a^\flat \in
QH_{n+k}(M)$ be its dual. Suppose that
$$
a^\flat = \sum_j a_j q^{A_j}
$$
with $ 0 \neq a_j \in H_{n+k + 2c_1(A_j)}(M)$ where $\lambda_j =
-\omega(A_j)$ are arranged as in (5.9). Denote by $\gamma$ a
Novikov cycle of $\e f$ with $[\gamma] = a^\flat \in HF_{k}(\e f)
\cong QH_{n + k}(M)$ and define the `gap'
$$
c(a) := \lambda_1 - \lambda_2.
$$
Then we have
$$
v(a) - {1 \over 2} c(a) \leq \inf_{\gamma}\{ \lambda_{\e
f}(\gamma) \mid [\gamma] = a^\flat\} \leq v(a) + {1 \over 2} c(a)
\tag 5.11
$$
for any sufficiently small $\e > 0$ and in particular, (5.8) holds.
We also have
$$
\lim_{\e \to 0}\inf_{\gamma}\{ \lambda_{\e f}(\gamma) \mid
[\gamma] = a^\flat\} = v(a) \tag 5.12
$$
and so the limit is independent of the choice of Morse functions $f$.
\endproclaim

\demo{Proof} We represent $a^\flat$ by a Novikov cycle
$$
\gamma = \sum_A\gamma_{A} q^{A}, \quad \gamma_j \in CM_*(-\e f)
$$
of $\e f$. It follows from (5.3) that if $A \in \Gamma(a)$,
all the coefficient Morse chains in
this sum must be cycles ,and if $A \not\in \Gamma(a)$, the
corresponding coefficient cycle must be a boundary. Therefore we can
decompose $\gamma$ as
$$
\gamma = \gamma_{\Gamma(a)} + \gamma_{\widetilde{ \Gamma(a)}} \tag
5.13
$$
where
$$
\align \gamma_{\Gamma(a)}& := \sum_{A \in \Gamma(a)} \gamma_A
q^{A}\\
\gamma_{\widetilde{\Gamma(a)}}& :=\sum_{B \not \in \Gamma(a)} \gamma_B
q^{B}
\endalign
$$
and we have
$$
\gamma_{\widetilde{\Gamma(a)}} = \part_{\e f}(\nu)
$$
for some Floer chain $\nu$ of $\e f$. Since the summands in
$\gamma_{\widetilde{ \Gamma(a)}}$ cannot cancel those in
$\gamma_{\Gamma(a)}$, we have
$$
\lambda_{\e f}(\gamma) \geq \lambda_{\e f}(\gamma_{\Gamma(a)}) =
\lambda_{\e f}\Big(\sum_{A \in \Gamma(a)} \gamma_A q^{A}\Big).
$$
Therefore by removing the exact term $\part_{\e f}(\gamma)$ when
we take the infimum over the cycles $\gamma$ with $[\gamma]
=a^\flat$ for the definition of $\rho(\e f;a)$, we may always
assume that $\gamma$ has the form
$$
\gamma = \sum_j \gamma_j q^{A_j}
$$
with $A_j \in \Gamma(a)$. Then again by (5.3), we have
$$
[\gamma_j] = a_j \in H_*(M).
$$
Furthermore we note that we have
$$
- \omega(A_j) - \max (\e f) \leq \lambda_{\e f}(\gamma_j q^{A_j})
\leq - \omega(A_j) - \min(\e f).
$$
Therefore if we choose $\e > 0$ so small that
$$
\e (\max f - \min f) \leq c(a) = \lambda_1 - \lambda_2,
$$
then we have
$$
\lambda_{\e f}(\gamma_1 q^{A_1}) \geq \lambda_{\e f}(\gamma_j
q^{A_j})
$$
for all $j = 1, \, 2, \cdots$ and so
$$
\lambda_{\e f}(\gamma) = \lambda_{\e f}(\gamma_1 q^{A_1}).
$$
Combining these, we derive
$$
- \omega(A_1) - \e \max f \leq \lambda_{\e f}(\gamma) \leq -
\omega(A_1) + \e \max f. \tag 5.14
$$
(5.11) follows from (5.14) if we choose $\e$ so that $\e (\max f -
\min f) < {c(a) \over 2}$. (5.12) also immediately follows from
(5.14). \qed\enddemo

\n {\it Step 2: The finiteness for general $H$.} Now we consider
generic nondegenerate $H$'s.
We fix $f$ be any Morse function and and $\e
> 0$ as in Lemma 5.4. Let $\alpha \in CF_*(H)$ be a Floer cycle of
$H$ with $[\alpha] = a^\flat$, and $\HH = \HH_{lin}$ the linear
homotopy
$$
\HH_{lin}: s \mapsto (1-s)(\e f) + s H.
$$
Applying (4.12) to the `inverse' linear homotopy
$$
\HH_{lin}^{-1}: s \mapsto (1-s)H + s (\e f).
$$
we obtain the inequality
$$
\lambda_{\e f}(h_{\HH_{lin}^{-1}}(\alpha)) \leq \lambda_H(\alpha)
+ \int_0^1 - \min (\e f - H)\, dt: \tag 5.15
$$
More precisely, it follows from the definition of the chain map
$h_{\HH_{lin}^{-1}}$ that for any generator $[z',w']$ of the cycle
$h_{\HH_{lin}^{-1}}(\alpha)$ of $\e f$, there is a generator
$[z,w]$ of the cycle $\alpha$ such that the equation (4.6) has a
solution. Then we derive from (4.11)
$$
\align \AA_{\e f}([z',w']) & \leq  \AA_H([z,w]) + \int_0^1
-\min_{x \in M} (\e f - H_t) \, dt\\
& \leq \lambda_H(\alpha) + \int_0^1 -\min_{x \in M} (\e f - H_t)
\, dt. \endalign
$$
Since this holds for any generator $[z',w']$ of
$\HH_{lin}^{-1}(\alpha)$, we obtain
$$
\lambda_{\e f}(\HH_{lin}^{-1}(\alpha)) \leq \lambda_H(\alpha) +
\int_0^1 -\min_{x \in M} (\e f - H_t) \, dt. \tag 5.16
$$
On the other hand, it follows from $[\HH_{lin}^{-1}(\alpha)] =
a^\flat$, that we have
$$
\rho(a;\e f) \leq \lambda_{\e f}(\HH_{lin}^{-1}(\alpha)).
$$
Combining this with (5.16), we derive
$$
\lambda_H(\alpha) \geq \rho(a;\e f) - \Big(\int_0^1 -\min_{x \in
M} (\e f - H_t) \, dt\Big). \tag 5.17
$$
Since this holds for any cycle $\alpha$ of $H$ with $[\alpha] =
a^\flat$, by taking the infimum over all such $\alpha$ in (5.17),
we have finally obtained
$$
\rho(H;a) \geq \rho(a;\e f) - \Big(\int_0^1 -\min_{x \in M} (\e f
- H_t) \, dt\Big). \tag 5.18
$$
Since Lemma 5.4 shows that $\rho(a;\e f) > -\infty$, this in
particular implies that $\rho(H;a) > -\infty$ and so $\rho(H;a)$
is finite.
\medskip

\n {\it Step 3: Proof of (5.4).} Finally we prove the inequality
(5.4). For this purpose, we consider general generic nondegenerate
pairs $H,\, K$. Let $\delta
> 0$ be any given number. We choose a cycle $\alpha$ of
$H$ respectively so that $[\alpha]= a^\flat$ and
$$
\lambda_H(\alpha) \leq \rho(H;a) + \delta \tag 5.19
$$
We would like to emphasize that {\it this is possible,  because we
have already shown that $\rho(H;a) > -\infty$}.

By considering the linear homotopy $h^{lin}_{HK}$ from $H$ to $K$,
we derive
$$
\lambda_K(h_{HK}^{lin}(\alpha)) \leq \lambda_H(\alpha) + \int
-\min_x(K_t - H_t) \, dt. \tag 5.20
$$
On the other hand (5.19) implies
$$
\aligned
\lambda_H(\alpha) & + \int -\min_x(K_t - H_t) \, dt \\
& \leq \rho(H;a) + \delta + \int -\min_x(K_t - H_t) \, dt
\endaligned. \tag 5.21
$$
Since $[h_{HK}^{lin}(\alpha)] = a^\flat$, we have
$$
\lambda_K(h_{HK}^{lin}(\alpha)) \geq \rho(K;a) \tag 5.22
$$
by the definition of $\rho(K;a)$. Combining (5.20)-(5.22), we have
derived
$$
\rho(K;a) - \rho(H;a) \leq \delta + \int_0^1 - \min_x(K_t-H_t)\, dt.
$$
Since this holds for arbitrary $\delta$, we have derived
$$
\rho(K;a) - \rho(H;a) \leq \int_0^1 - \min_x(K_t-H_t)\, dt.
$$
By changing the role of $H$ and $K$, we also derive
$$
\rho(H;a) - \rho(K;a) \leq  \int_0^1 - \min_x(H_t-K_t)\, dt
= \int _0^1 \max_x (K_t-H_t)\, dt
$$
Hence, we have the inequality
$$
\int_0^1 - \max_x(K_t-H_t)\, dt  \leq \rho(K;a) -
\rho(H;a)
\leq \int_0^1 - \min_x(K_t-H_t)\, dt
$$
which is precisely (5.4).  Obviously the inequality (5.4),
enables us to
extend the definition of $\rho$ by continuity to arbitrary
$C^0$-Hamiltonians. This finishes the proof of Theorem 5.3.
\qed\enddemo

\head{\bf \S 6. Basic properties of the spectral invariants}
\endhead

In this section, we will prove all the remaining properties stated
in Theorem I in the introduction. We first re-state the main
axioms of the spectral invariants.

\proclaim{Theorem 6.1} Let $H,\, F$ be arbitrary smooth
Hamiltonian functions, and $a \neq 0 \in QH^*(M)$ and let
$$
\rho: C^\infty_m([0,1] \times M) \times QH^*(M) \to \R
$$
be as defined in \S 5.  Then $\rho$ satisfies the following
properties: \roster

\item {\bf (Projective invariance)} $\rho(H;\lambda a) =\rho(H;a)$
for any $0 \neq \lambda \in \Q$ \item {\bf (Normalization)} For $a
= \sum_{ A \in \Gamma} a_A \otimes q^A$, $\rho(\underline 0;a) =
v(a), \text{the valuation of }\, a$. \item {\bf (Symplectic
invariance)} $\rho(\eta^*H;a) = \rho(\widetilde H;a)$ for any
symplectic diffeomorphism $\eta$. \item {\bf (Triangle
inequality)} $\rho(H\# F; a\cdot b) \leq \rho(H;a) + \rho(F;b)$
\item {\bf ($C^0$-Continuity)} $|\rho(H;a) - \rho(F;a)| \leq \|H -
F \|$ and in particular $\rho(\cdot, a)$ is continuous with
respect to the $C^0$-topology of Hamiltonian functions.
\endroster
\endproclaim

We have already proven the properties of normalization and $C^0$
continuity in the course of proving the linking property of the
Novikov Floer cycles in \S 5.  The remaining parts of the proofs
deal with the {\it symplectic invariance}
and the {\it triangle inequality}.
\medskip

\n {\it 6.1. Proof of symplectic invariance.}
\smallskip

We consider the symplectic conjugation
$$
\phi \mapsto \eta^{-1}\phi \eta; \quad \HH am(M,\omega) \to \HH
am(M,\omega)
$$
for any symplectic diffeomorphism $\eta: (M,\omega) \to
(M,\omega)$. Recall that the pull-back function $\eta^*H$ given
by
$$
\eta^*H(t,x)=H(t, \eta(x)) \tag 6.1
$$
generates the conjugation $\eta^{-1}\phi\eta$ when $H \mapsto
\phi$.

We summarize the basic facts on this conjugation relevant to the
filtered Floer homology here:

\roster \item when $H \mapsto \phi$, $\eta^*H \mapsto \eta \phi
\eta^{-1}$, \item if $H$ is nondegenerate, $\eta^*H$ is also
nondegenerate, \item if $(J,H)$ is regular in the Floer theoretic
sense, then so is $(\eta^*J, \eta^*H)$, \item there exists natural
bijection $\eta_*: \Omega_0(M) \to \Omega_0(M)$ defined by
$$
\eta_*([z,w]) = ([\eta\circ z, \eta\circ w])
$$
under which we have the identity
$$
\AA_H([z,w]) = \AA_{\eta^*H}(\eta_*[z,w]). \tag 6.2
$$
\item the $L^2$-gradients of the corresponding action functionals
satisfy
$$
\eta_*(\hbox{grad}_J\AA_H)([z,w]) =
\hbox{grad}_{\eta^*J}(\AA_{\eta^*H})(\eta_*([z,w])) \tag 6.3
$$
\item if $u: \R \times S^1 \to M$ is a solution of perturbed
Cauchy-Riemann equation for $(J,H)$, then $\eta_*u = \eta\circ u$
is a solution for the pair $(\eta^*J, \eta^*H)$. In addition, all
the Fredholm properties of $(J,H,u)$ and $(\eta^*J,
\eta^*H,\eta_*u)$ are the same.
\endroster
These facts imply that the conjugation by $\eta$ induces the
canonical filtration preserving chain isomorphism
$$
\eta_*: (CF_*^\lambda(H), \part_{(H,J)}) \to (CF_*^\lambda
(\eta^*H), \part_{(\eta^*H,\eta^*J)})
$$
for any $\lambda \in \R \setminus \hbox{Spec}(H) = \R \setminus
\hbox{Spec}(\eta^*H)$. In particular it induces a filtration
preserving isomorphism
$$
\eta_*: HF_*^\lambda(H, J) \to HF_*^\lambda (\eta^*H,\eta^*J).
$$
in homology. The symplectic invariance is then an immediate
consequence of our construction of $\rho(H;a)$.

\medskip
\n{\it 6.2. Proof of the triangle inequality}
\smallskip

To start with the proof of the triangle inequality,  we need to
recall the definition of the ``pants product''
$$
HF_*(H) \otimes HF_*(F) \to HF_*(H \# F). \tag 6.4
$$
We also need to straighten out the grading problem of the pants
product.

For the purpose of studying the effect on the filtration under the
product, we need to define this product in the chain level in an
optimal way as in [Oh2], [Sc] and  [En1].  For this purpose, we
will mostly follow the description provided by Entov [En1] with
few notational changes and different convention on the grading. As
pointed out before, our grading convention satisfies (4.17) under
the pants product. Except the grading convention, the conventions
in [En1,2] on the definition of Hamiltonian vector field and the
action functional coincide with our conventions in [Oh1-3,5] and
here.

Let $\Sigma$ be the compact Riemann surface of genus 0
with three punctures. We fix a holomorphic identification of a
neighborhood of each puncture with either $[0, \infty) \times S^1
$ or $(-\infty, 0] \times S^1$ with the standard complex structure
on the cylinder. We call punctures of the first type {\it
negative} and the second type {\it positive}. In terms of the
``pair-of-pants'' $\Sigma \setminus \cup_i D_i$, the positive
puncture corresponds to the {\it outgoing ends} and the negative
corresponds to the {\it incoming ends}. We denote the
neighborhoods of the three punctures by $D_i$, $i = 1, 2, 3$ and
the identification by
$$
\varphi^+_i: D_i \to (-\infty, 0] \times S^1
$$
for positive punctures and
$$
\varphi^-_3: D_3 \to [0, \infty) \times S^1
$$
for negative punctures.
We denote by $(\tau,t)$ the standard cylindrical coordinates on
the cylinders.

We fix a cut-off function $\rho^+: (-\infty,0] \to [0,1]$
defined by
$$
\rho = \cases 1 \quad & \tau \leq -2 \\
0 \quad & \tau \geq -1
\endcases
$$
and $\rho^-: [0,\infty) \to [0,1]$ by $\rho^-(\tau) =
\rho^+(-\tau)$. We will just denote by $\rho$ these cut-off
functions for both cases when there is no danger of confusion.

We now consider the (topologically) trivial bundle $P \to \Sigma$ with fiber
isomorphic to $(M,\omega)$ and fix a trivialization
$$
\Phi_i: P_i:= P|_{D_i} \to D_i \times M
$$
on each $D_i$.
On each $P_i$, we consider the closed two form of the type
$$
\omega_{P_i}:= \Phi_i^*(\omega + d(\rho H_t dt)) \tag 6.5
$$
for a time periodic Hamiltonian $H: [0,1] \times M \to \R$. The
following is an important lemma whose proof we omit (see [En1]).

\proclaim{Lemma 6.2} Consider three normalized Hamiltonians $H_i$,
$i = 1, 2,3$. Then there exists a closed 2-form $\omega_P$ such
that \roster \item $\omega_P|_{P_i} = \omega_{P_i}$ \item
$\omega_P$ restricts to $\omega$ in each fiber \item
$\omega_P^{n+1} = 0$
\endroster
\endproclaim
Such $\omega_P$ induces a canonical symplectic connection $\nabla
= \nabla_{\omega_P}$ [GLS], [En1]. In addition it also fixes a natural
deformation class of symplectic forms on $P$ obtained by those
$$
\Omega_{P,\lambda} := \omega_P + \lambda \omega_\Sigma
$$
where $\omega_\Sigma$ is an area form and $\lambda > 0$ is a
sufficiently large constant. We will always normalize
$\omega_\Sigma$ so that $\int_\Sigma \omega_\Sigma = 1$.

Next let $\widetilde J$ be an almost complex structure on $P$ such
that
\roster
\item $\widetilde J$ is $\omega_P$-compatible on each fiber and so preserves
the vertical tangent space
\item the projection $\pi: P \to \Sigma$ is pseudo-holomorphic, i.e,
$d\pi \circ \widetilde J = j \circ d\pi$.
\endroster
When we are given three $t$-periodic Hamiltonian $H =
(H_1,H_2,H_3)$, we say that $\widetilde J$ is $(H,J)$-compatible,
if $\widetilde J$ additionally satisfies

(3) For each $i$, $(\Phi_i)_*\widetilde J = j\oplus J_{H_i}$ where
$$
J_{H_i}(\tau,t,x) = (\phi^t_{H_i})^*J
$$
for some $t$-periodic family of almost complex structure $J = \{
J_t\}_{0 \leq t \leq 1}$ on $M$ over a disc $D_i' \subset D_i$ in
terms of the cylindrical coordinates. Here $D_i' =
\varphi_i^{-1}((-\infty, -K_i] \times S^1), i = 1,\, 2$ and
$\varphi^{-1}_3([K_3, \infty) \times S^1)$ for some $K_i > 0$.
Later we will particularly consider the case where $J$ is in the
special form adapted to the Hamiltonian $H$. See (6.23).
\medskip

The condition (3) implies that the $\widetilde J$-holomorphic sections
$v$ over $D_i'$ are precisely the solutions of the equation
$$
{\part u \over \part \tau} + J_t\Big({\part u \over \part t} -
X_{H_i}(u)\Big) = 0 \tag 6.6
$$
if we write $v(\tau,t) = (\tau,t, u(\tau,t))$ in the trivialization
with respect to the cylindrical coordinates $(\tau,t)$ on $D_i'$
induced by $\phi_i^\pm$ above.

Now we are ready to define the moduli space which will be relevant
to the definition of the pants product that we need to use. To
simplify the notations, we denote
$$
\widehat z =[z,w]
$$
in general and $\widehat z = (\widehat z_1, \widehat z_2, \widehat
z_3)$ where $\widehat z_i =[z_i,w_i] \in \text{Crit}\AA_{H_i}$ for
$i =1,2, 3$.

\definition{Definition 6.3} Consider the Hamiltonians
$H =\{H_i\}_{1\leq i \leq 3}$ with $H_3 = H_1 \# H_2$, and let
$\widetilde J$ be a $H$-compatible almost complex structure. We
denote by $\MM(H, \widetilde J; \widehat z)$ the space of all
$\widetilde J$-holomorphic sections $u: \Sigma \to P$ that satisfy
\roster \item The maps $u_i:= u \circ (\phi_i^{-1}): (-\infty,
K_i] \times S^1 \to M$ which are solutions of (6.6), satisfy
$$
\lim_{\tau \to -\infty}u_i(\tau,\cdot) = z_i, \quad i = 1,2
$$
and similarly for $i=3$ changing $-\infty$ to $+\infty$.
\item
The closed surface obtained by capping off $pr_M\circ u(\Sigma)$
with the discs $w_i$ taken with the same orientation for $i = 1,2$
and the opposite one for $i =3$ represents zero (mod) $\Z$-torsion
elements.
\endroster
\enddefinition

Note that $\MM(H, \widetilde J; \widehat z)$ depends only on the
equivalence class of $\widehat z$'s: we say that $\widehat z' \sim
\widehat z$ if they satisfy
$$
z_i' = z_i, \quad w_i' = w_i \# A_i
$$
for $A_i \in \pi_2(M)$ and $\sum_{i=1}^3 A_i$ represents zero
(mod) $\Z$-torsion elements. The (virtual) dimension of $\MM(H,
\widetilde J; \widehat z)$ is given by
$$
\aligned \dim  \MM(H, \widetilde J; \widehat z) & = 2n -
(\mu_{H_1}(z_1) + n) - (\mu_{H_2}(z_2) +n) - (-\mu_{H_3}(z_3) + n)\\
& = n + (\mu_{H_3}(z_3) - \mu_{H_1}(z_1) - \mu_{H_2}(z_2)).
\endaligned
\tag 6.7
$$
Note that when $\dim \MM(H,\widetilde J;\widehat z) = 0$, we have
$$
n= -\mu_{H_3}(\widehat z_3)+\mu_{H_1}(\widehat z_1)+\mu_{H_2}(\widehat z_2)
$$
which is equivalent to
$$
\mu_{H_3}(\widehat z_3) =(\mu_{H_1}(\widehat
z_1)+\mu_{H_2}(\widehat z_2)) -n \tag 6.8
$$
which provides the degree of the
pants product (4.17) in our convention of the grading of
the Floer complex we adopt in the
present paper.  Now the pair-of-pants product $*$ for the
chains is defined by
$$
\widehat z_1 * \widehat z_2 =  \sum_{\widehat z_3} \#(\MM(H,
\widetilde J;\widehat z)) \widehat z_3 \tag 6.9
$$
for the generators $\widehat z_i$ and then by linearly extending over
the chains in $CF_*(H_1) \otimes CF_*(H_2)$. Our grading convention makes
this product is of degree zero. Now with this preparation, we are ready
to prove the triangle inequality.

\demo{Proof of the triangle inequality} Let $\alpha \in CF_*(H)$
and $\beta \in CF_*(F)$ be Floer cycles with $[\alpha] = [\beta] =
a^\flat$ and consider their pants product cycle $\alpha * \beta:
=\gamma \in CF_*(H\# F)$. Then we have
$$
[\alpha * \beta] = (a\cdot b)^\flat
$$
and so
$$
\rho(H\# F;a\cdot b) \leq \lambda_{H\# F}(\alpha * \beta). \tag
6.10
$$
Let $\delta > 0$ be any given number and choose $\alpha\in
CF_*(H)$ and $\beta \in CF_*(F)$ so that
$$
\aligned
\lambda_H(\alpha) & \leq \rho(H;a) + \frac{\delta}{2} \\
\lambda_H(\beta) & \leq \rho(F;b) + \frac{\delta}{2}.
\endaligned
\tag 6.11
$$
Then we have the expressions
$$
\alpha = \sum_i a_i [z_i,w_i] \, \text{ with }\,
\AA_{H}([z_i,w_i]) \leq \rho(H;a) + \frac{\delta}{2}
$$
and
$$
\beta = \sum_j a_j [z_j,w_j] \, \text{ with } \,
\AA_{H}([z_j,w_j]) \leq \rho(H;b) + \frac{\delta}{2}.
$$
Now using the pants product (6.9), we would like to estimate the
level of the chain $\alpha \# \beta \in CF_*(H\# F)$. The
following is a crucial lemma whose proof we omit but refer to
[Sect. 4.1, Sc] or [Sect. 5, En1].

\proclaim{Lemma 6.4} Suppose that $\MM(H,\widetilde J;\widehat z)$
is non-empty. Then we have the identity
$$
\int v^*\omega_P = - \AA_{H_1\# H_2}([z_3,w_3]) +
\AA_{H_1}([z_1,w_1]) + \AA_{H_2}([z_2,w_2]) \tag 6.12
$$
for any  $\in \MM(H,\widetilde J;\widehat z)$
\endproclaim

Now since $\widetilde J$-holomorphic and $\widetilde J$
is $\Omega_{P,\lambda}$-compatible, we have
$$
0 \leq \int v^*\Omega_{P,\lambda} = \int v^*\omega_P + \lambda
\int v^*\omega_\Sigma = \int v^*\omega_P + \lambda.
$$

\proclaim{Lemma 6.5 [Theorem 3.6.1 \& 3.7.4, En1]} Let $H_i$'s be
as in Lemma 6.2. Then for any given $\delta > 0$, we can choose a
closed 2-form $\omega_P$ so that $\Omega_{P,\lambda} = \omega_P +
\lambda \omega_\Sigma$ becomes a symplectic form for all $\lambda
\geq \delta$. In other words, the size $size(H)$ (see [Definition
3.1, En1]) is $\infty$.
\endproclaim

We recall that from the definition of $*$ that for any $[z_3,w_3]
\in \alpha
* \beta$ there exist $[z_1, w_1] \in \alpha$ and $[z_2,w_2] \in
\beta$ such that $\MM(\widetilde J, H;\widehat z)$ is non-empty
with the asymptotic condition
$$
\widehat z = ([z_1,w_1], [z_2,w_2];[z_3,w_3]).
$$
Applying this and the above two lemmata to $H$ and $F$ for
$\lambda$ arbitrarily close to 0, and also applying (6.10) and
(6.11),  we immediately derive
$$
\align \AA_{H\# F}([z_3,w_3]) & \leq \AA_H([z_1,w_1]) +
\AA_F([z_2,w_2]) + \delta \\
& \leq \lambda_H(\alpha) + \lambda_F(\beta) + \delta \\
& \leq \rho(H;a) + \rho(F;b) + 2\delta \tag 6.13
\endalign
$$
for any $[z_3,w_3] \in \alpha * \beta$. Combining (6.10), (6.11)
and (6.13),  we derive
$$
\rho(H\# F;a\cdot b) \leq \rho(H;a) + \rho(F;b) + 2\delta
$$
Since this holds for any $\delta$, we have proven
$$
\rho(H\# F;a\cdot b) \leq \rho(H;a) + \rho(F;b).
$$
The triangle inequality mentioned in Theorem 6.1 immediately
follows from the definition $\rho(\widetilde\phi;a) = \rho(H;a)$
in Theorem 5.5. \qed\enddemo

\head{\bf \S 7. The rational case; proof of the spectrality}
\endhead

In this section, we will prove the spectrality for the rational
sympelctic manifolds: we recall that a symplectic manifold
$(M,\omega)$ {\it rational} if the period group $\Gamma_\omega$ is
discrete. We will further study the spectrality property on
general symplectic manifolds elsewhere, which turns out to be much
more nontrivial to prove.

\proclaim{Theorem 7.1} Suppose that $(M,\omega)$ is rational. Then
for any smooth one-periodic Hamiltonian function $H: S^1 \times M
\to \R$, we have
$$
\rho(H;a) \in \text{\rm Spec}(H)
$$
for each given quantum cohomology class $0 \neq a \in QH^*(M)$.
\endproclaim

\demo{Proof }  We need to show that the mini-max
value $\rho(H;a)$ is a critical value, i.e., that there exists
$[z,w] \in \widetilde \Omega_0(M)$ such that
$$
\aligned
& \AA_H([z,w]) = \rho(H;a) \\
& d\AA_H([z,w]) = 0, \quad \text{i.e., } \quad \dot z = X_H(z).
\endaligned
\tag 7.1
$$
We have already shown the finiteness of the value  $\rho(H;a)$ in section 5.
If $H$ is nondegenerate, we just use the fixed Hamiltonian.
If $H$ is not nondegenerate, we approximate $H$ by a sequence of nondegenerate
Hamiltonians $H_j$ in the $C^2$ topology.  Let $[z_j,w_j] \in
\text{Crit}\AA_{H_j}$ be the peak of a Floer cycle $\alpha_j
\in CF_*(H_j)$,  such that
$$
\lim_{j \to \infty} \AA_{H_j}([z_j,w_j]) = \rho(H;a). \tag 7.2
$$
Such a sequence can be chosen from the definition of $\rho(\cdot; a)$ and
the finiteness thereof.

Since $M$ is compact and $H_j \to H$ in the $C^2$ topology,  and
$\dot z_j = X_{H_j}(z_j)$ for all $j$, it follows from the
standard boot-strap argument that $z_j$ has a subsequence, which
we still denote by $z_j$, converging to $z_\infty$ which solves
$\dot z = X_H(z)$. Now we show that $[z_j,w_j]$ themselves are
pre-compact on $\widetilde \Omega_0(M)$.  Since we fix the quantum
cohomolgy class $0 \neq a \in QH^*(M)$ (or more specifically since we
fix its degree) and the Floer cycle satisfies $[\alpha_j] = a$, we have
$$
\mu_{H_j}([z_j,w_j]) = \mu_{H_i}([z_i,w_i]). \tag 7.3
$$

\proclaim{Lemma 7.2} When $(M,\omega)$ is rational,
$\text{Crit}(\AA_K) \subset \widetilde \Omega_0(M)$ is a closed
subset of $\R$ for any smooth Hamiltonian $K$, and is locally
compact in the subspace topology  of  the covering space
$$
\pi: \widetilde \Omega_0(M) \to \Omega_0(M).
$$
\endproclaim
\demo{Proof} First note that when $(M,\omega)$ is rational, the covering group
$\Gamma_\omega$ of $\pi$ above is discrete.   Together with the fact that the set
of solutions of $\dot z = X_K(z)$ is compact (on compact M),  it
follows that
$$
\text{Crit}(\AA_K) =  \{[z,w] \in \widetilde\Omega_0(M)
 \mid \dot z = X_K(z) \}
$$
is a closed subset which is also locally compact. \qed\enddemo

Now consider the bounding discs of $z_\infty$
$$
w_j' = w_j \# u^{can}_j
$$
for  all sufficiently large $j$, where $u^{can}$ is the
homotopically unique thin cylinder between $z_j$ and $z_\infty$:
more precisely, $u_j^{can}$ is given by the formula
$$
u_j^{can}(s,t) = \exp_{z_j}(s \xi_j(t)), \quad \xi_j(t) =
(\exp_{z_j(t)})^{-1} (z_\infty(t)) \tag 7.4
$$
where $\exp$ is the exponential map with respect to a fixed metric
$g_{J_{ref}} = \omega(\cdot, J_{ref}\cdot)$ for a fixed compatible
almost complex structure. We note that as $j \to \infty$ the
geometric area of $u_j^{can}$ converges to $0$.

We compute the action of the critical points $[z_\infty, w_j'] \in
\text{Crit}\AA_H$,
$$
\align \AA_H([z_\infty, w_j'])  & = - \int_{w_j'} \omega -
\int_0^1 H(t, z_\infty(t)) \, dt
\tag 7.5 \\
& = - \int_{w_j} \omega - \int_{u^{can}_j} - \int_0^1 H(t, z_\infty(t)) \, dt \\
& = \Big(- \int_{w_j} \omega  - \int_0^1 H_j(t, z_j(t)) \, dt\Big) \\
& \qquad - \Big(\int_0^1 H(t, z_\infty(t))
- \int_0^1 H_j(t, z_j(t)) \Big) - \int_{u^{can}_j} \omega. \\
& = \AA_{H_j}([z_j, w_j]) - \Big(\int_0^1 H(t, z_\infty(t))
- \int_0^1 H_j(t, z_j(t)) \Big) - \int_{u^{can}_j} \omega. \\
\endalign
$$
From the explicit expression (7.4), it follows that
$$
\lim_{j \to \infty} \int_{u^{can}_j}\omega = 0
\tag 7.6
$$
since the geometric area of $u_j^{can}$ converges to zero and we
have $\text{Area}(u_j^{can}) \geq |\int_{u^{can}_j}\omega|$. Since
$z_j$ converges to $z_\infty$ uniformly and $H_j \to H$, we have
$$
- \Big(\int_0^1 H(t, z_\infty(t)) - \int_0^1 H(t, z_j(t)) \Big)
\to 0. \tag 7.7
$$
Therefore combining (7.2), (7.6) and (7.7), we derive
$$
\lim_{j \to \infty} \AA_H([z_\infty, w_j'])  =  \rho(H;a).
$$
In particular $\AA_H([z_\infty, w_j'])$ is a Cauchy sequence,  which implies
$$
\Big|\int_{w_j'} \omega - \int_{w_i'}\omega \Big| =
\Big|\AA_H([z_\infty, w_j']) - \AA_H([z_\infty, w_i'])\Big| \to 0
$$
i.e.,
$$
\int_{w_j' \# \overline w_i'} \omega \to 0.
$$
Since $\Gamma_\omega$ is discrete and $\int_{w_j' \# \overline
w_i'} \omega\in \Gamma_\omega$, this indeed implies that
$$
\int_{w_j' \# \overline w_i'} \omega  = 0 \tag 7.8
$$
for all sufficiently large $i, \, j$. Since the set
$\Big\{\int_{w_j'} \omega\Big\}_{j \in\Z_+}$ is bounded, these
imply that the sequence $\int_{w_j'} \omega$ eventually
stabilizes. Going back to (7.5), we have proven that the actions
$$
\AA_H([z_\infty, w_j'])
$$
themselves stabilize and so we have
$$
\AA_H([z_\infty,w_N']) = \lim_{j \to \infty}\AA_H([z_\infty,
w_j']) = \rho(H;a)
$$
for a fixed sufficiently large $N \in \Z_+$. This proves that
$\rho(H;a)$ is indeed a critical value of $\AA_H$ at the critical
point $[z_\infty,w_N']$. This finishes the proof. \qed\enddemo

We now state the following theorem.

\proclaim{Theorem 7.3} Let $(M,\omega)$ be rational and
$C_m^\infty(M\times [0,1],\R)$ be the set of normalized
$C^\infty$-Hamiltonians on $M$. We denote by
$\rho_a:C_m^\infty(M\times [0,1],\R) \to \R$ the extended
continuous function defined by $\rho_a(H) = \rho(H; a)$.\roster

\item The image of $\rho_a$ depends only on the homotopy class
$\widetilde \phi= [\phi,H]$. Hence $\rho_a$ pushes down to a
well-defined function
$$
\rho: \widetilde{\HH am}(M,\omega) \times QH^*(M) \to \R; \quad
\rho(\widetilde \phi;a) : = \rho(H;a) \tag 7.9
$$
for any $H$ with $\widetilde \phi = [\phi,H]$.

\item We have the formula
$$
\rho(H;a) = \inf_{\lambda}\{ \lambda \mid a^\flat \in \text{Im
}(i_\lambda: HF_*^\lambda(H) \to HF_*(H)) \}. \tag 7.10
$$
\endroster
\endproclaim
\demo{Proof}  We have shown in Theorem 7.1 that $\rho(H;a)$ is
indeed a critical value of $\AA_H$, i.e., lies in
$\text{Spec}(H)$.  With this fact in our disposition, the
well-definedness of the definition (7.9), i.e., independence of
$H$ with $\widetilde \phi = [\phi,H]$ is an immediate consequence
of combination of the following results: \roster \item $H \mapsto
\rho(H;a)$ is continuous, \item $\text{Spec}(H)$ is of measure
zero subset of $\R$ (Lemma 2.2), \item $\text{Spec}(H) =
\text{Spec}(\widetilde \phi)$  depends only on its homotopy class
$[H]=\widetilde \phi$ and so fixed as long as $[H] = \widetilde
\phi$ (Theorem 2.3), \item  the only real-valued  continuous
function from a connect space (e.g., $[0,1]$) whose image has
measure zero in $\R$, is  a constant function.
\endroster

(7.10) is just a rephrasing of the definition of $\rho(H;a)$. This
finishes the proof of Theorem 7.3. \qed\enddemo

One more important property concerns the effect of $\rho$ under
the action of $\pi_0(\widetilde G)$. We first explain how
$\pi_0(\widetilde G)$ acts on $\widetilde{\HH am}(M,\omega) \times
QH^*(M)$ following (and adapting into cohomological version)
Seidel's description of the action on $QH_*(M)$. According to
[Se], each element $[h,\widetilde h] \in \pi_0(\widetilde G)$ acts
on $QH_*(M)$ by the quantum product of an even element
$\Psi([h,\widetilde h])$ on $QH_*(M)$. We take the adjoint action
of it on $a \in QH^*(M)$ and denote it by $\widetilde h^*a$. More
precisely, $\widetilde h^*a$ is defined by the identity
$$
\langle \widetilde h^* a, \beta \rangle =\langle a,
\Psi([h,\widetilde h])\cdot \beta \rangle \tag 7.11
$$
with respect the non-degenerate pairing $\langle \cdot, \cdot
\rangle$ between $QH^*(M)$ and $QH_*(M)$.

\proclaim{Theorem 7.4} {\bf (Monodromy shift)}  Let $\pi_0(\widetilde G)$ act on
$\widetilde {\HH am}(M,\omega) \times QH^*(M)$ as above, i.e,
$$
[h,\widetilde h] \cdot (\widetilde \phi,a) = (h\cdot \widetilde
\phi, \widetilde h^*a) \tag 7.12
$$
Then we have
$$
\rho([h,\widetilde h]\cdot
(\widetilde \phi;a))
 = \rho(\widetilde\phi;a) + I_\omega([h,\widetilde h]).
$$
\endproclaim
\demo{Proof} This is immediate from the construction of
$\Psi([h,\widetilde h])$ in [Se]. Indeed, the map
$$
[h, \widetilde h]_*: CF_*(F) \mapsto CF_*(H\# F) \tag 7.13
$$
is induced by the map
$$
[z,w] \mapsto \widetilde h([z,w])
$$
and we have
$$
\AA_{H \# F}(\widetilde h ([z,w])) = \AA_F([z,w]) +
I_\omega([h,\widetilde h])
$$
by (2.5). Furthermore the map (7.12) is a chain isomorphism whose
inverse is given by $([h,\widetilde h]^{-1})_*$. This immediately
implies the theorem from the construction of $\rho$. \qed\enddemo

\definition{Remark 7.5} Strictly speaking, $\widetilde h^*a$ may
not lie in the standard quantum cohomology $QH^*(M)$ because it is
defined as the linear functional on the complex $CQ_*(M)$ that is
dual to the Seidel element $\Psi([h,\widetilde h]) \in CQ_*(M)$
under the canonical pairing between $CQ_*(M)$ and $CQ^*(M)$. A
priori, the {\it bounded} linear functional
$$
\widetilde h^*a = \langle \Psi([h,\widetilde h]), \cdot \rangle
$$
may not lie in the image of $\sharp: QH_*(M) \to QH^*(M)$,
mentioned in section 3, in general. In that case, one
should consider $\widetilde h^*a$ as a {\it continuous} quantum
cohomology class in the sense of Appendix. We refer readers to
Appendix  for the explanation on how to extend the definition of
our spectral invariants to the continuous quantum cohomology
classes.
\enddefinition

Now we can define
$$
\rho: \widetilde{\HH am}(M,\omega) \times QH^*(M)
$$
by putting
$$
\rho(\widetilde \phi;a): = \rho(H;a)
$$
for any $H \mapsto \widetilde \phi$ with $[H] = \widetilde \phi$
when $\widetilde \phi$ is nondegenerate, and then extending to
arbitrary $\widetilde\phi$ by continuity.

Then by the spectrality of $\rho(\widetilde \phi;a)$ for each $a
\in QH^*(M)$, we have constructed continuous `sections' of the
action spectrum bundle
$$
\frak{Spec}(M,\omega) \to \widetilde{\HH am}(M,\omega)
$$
We define the {\it essential spectrum} of $\widetilde \phi$ by
$$
\align \text{spec}(\widetilde \phi)
& : = \{\rho(\widetilde \phi;a) \mid 0 \neq a \in QH^*(M)\} \\
\text{spec}_k(\widetilde \phi) & : = \{\rho(\widetilde \phi;a)
\mid 0 \neq a \in QH^k(M)\}
\endalign
$$
and the bundle of essential spectra by
$$
\frak{spec}(M,\omega) = \bigcup_{\widetilde \phi \in
\widetilde{\HH am}(M,\omega)} \text{spec}(\widetilde \phi)
$$
and similarly for $\frak{spec}_k(M,\omega)$.

\head{\bf \S 8. Remarks on the transversality}
\endhead

Our construction of various maps in the Floer homology works as
they are in the previous section for the strongly semi-positive
case [Se], [En1] by the standard transversality argument. On the
other hand in the general case where constructions of operations
in the Floer homology theory requires the machinery of virtual
fundamental chains through multi-valued abstract perturbation
[FOn], [LT1], [Ru], we need to explain how this general machinery
can be incorporated in our construction. The full details will be provided
elsewhere.  We will use the
terminology `Kuranishi structure' adopted by Fukaya and Ono [FOn]
for the rest of the discussion.

One essential point in our proofs is that various numerical
estimates concerning the critical values of the action functional
and the levels of relevant Novikov cycles do {\it not} require
transversality of the solutions of the relevant pseudo-holomorphic
sections, but {\it depends only on the non-emptiness of the
moduli space}
$$
\MM(H, \widetilde J;\widehat z)
$$
which can be studied for {\it any}, not necessarily generic,
Hamiltonian $H$. Since we always have suitable a priori energy
bound which requires some necessary homotopy assumption on the
pseudo-holomorphic sections, we can compactify the corresponding
moduli space into a compact Hausdorff space, using a variation of
the notion of stable maps in the case of non-degenerate
Hamiltonians $H$. We denote this compactification again by
$$
\MM(H,\widetilde J;\widehat z).
$$
This space could be pathological in general.
But because we assume that the Hamiltonians $H$ are non-degenerate, i.e,
all the periodic orbits are non-degenerate, the moduli space
is not completely pathological but at least carries a Kuranishi
structure in the sense of Fukaya-Ono [FOn] for any $H$-compatible
$\widetilde J$.  This enables us to
apply the abstract multi-valued perturbation theory and
to perturb the compactified moduli space by a Kuranishi map
$\Xi$ so that the perturbed  moduli space
$$
\MM(H, \widetilde J;\widehat z, \Xi)
$$
is transversal in that the linearized equation of the perturbed
equation
$$
\overline\part_{\widetilde J}(v) + \Xi(v) = 0
$$
is surjective and so its solution set carries a smooth (orbifold)
structure. Furthermore the perturbation $\Xi$ can be chosen so
that as $\|\Xi\| \to 0$, the perturbed moduli space $\MM(H,
\widetilde J;\widehat z, \Xi)$ converges to $\MM(H, \widetilde J;
\widehat z)$ in a suitable sense (see [FOn] for the precise
description of this convergence).

Now the crucial point is that non-emptiness of the perturbed
moduli space will be guaranteed as long as certain topological
conditions are met. For example, the followings are the prototypes
that we have used in this paper: \roster
\item $h_{\HH_1}: CF_0(\e f) \to CF_0(H)$ is an isomorphism in
homology and so $[h_{\HH_1}(1^\flat)] \neq 0$. This is immediately
translated as an existence result of solutions of the
perturbed Cauchy-Riemann equation.
\item
The definition of the pants product $*$ and the identity
$$
[\alpha * \beta] =  (a\cdot b)^\flat
$$
in homology guarantee non-emptiness of the relevant perturbed
moduli space $\MM(H,\widetilde J;\widehat z, \Xi)$ for $\alpha \in
CF_*(H_1), \, \beta \in CF_*(H_2)$ with $[\alpha] = a^\flat$ and
$[\beta] = b^\flat$ respectively.
\endroster

Once we prove non-emptiness of $\MM(H,\widetilde J;\widehat z, \Xi)$
and an a priori energy bound for the
non-empty perturbed moduli space and {\it if the asymptotic conditions
$\widehat z$ are fixed}, we can study the convergence
of a sequence $v_j \in \MM(H, \widetilde J; \widehat z, \Xi_j)$
as $\Xi_j \to 0$ by the Gromov-Floer compactness theorem.
However a priori there are infinite
possibility of asymptotic conditions for the pseudo-holomorphic
sections that we are studying, because we typically impose
that the asymptotic limit lie in certain Novikov
cycles like
$$
\widehat z_1 \in \alpha, \, \widehat z_2 \in
\beta, \, \widehat z_3 \in \alpha*\beta
$$
Because the Novikov Floer cycles are generated by an infinite number of
critical points $[z,w]$ in general, one needs to control the
asymptotic behavior to carry out compactness argument. For this
purpose, we need to establish a {\it lower bound} for the actions
which will enable us to consider only finite possibilities for the
asymptotic conditions because of the finiteness condition in the
definition of Novikov chains. We would like to emphasize that
obtaining a lower bound is the heart of matter in the classical
mini-max theory of the {\it indefinite} action functional which
requires a linking property of semi-infinite cycles. On the other
hand, obtaining {\it upper bound} is usually an immediate
consequence of the identity like (4.10).

With such a lower bound for the actions, we may then assume, by
taking a subsequence if necessary, that the asymptotic conditions
are fixed when we take the limit and so we can safely apply the
Gromov-Floer compactness theorem to produce a (cusp)-limit lying
in the compactified moduli space $\MM(H, \widetilde J; \widehat
z)$. This will then justify all the statements and proofs in the
previous sections for the complete generality.

\vskip0.5truein

\head{\bf Appendix: Continuous quantum cohomology}
\endhead

In this appendix, we define the genuinely cohomological version of
the quantum cohomology and explain how we can extend the
definition of the spectral invariants to the classes in this
cohomolgical version.

 We call this
{\it continuous quantum cohomology} and denote by
$$
QH^*_{cont}(M).
$$
In this respect, we call the usual quantum cohomology ring
$QH^*(M) = H^*(M)\otimes \Lambda^\uparrow$ the {\it finite quantum
cohomology}. We call elements in $QH^*_{cont}(M)$ and $QH^*(M)$
continuous (resp. finite) quantum cohomology classes.

We first define the chain complex associated to $QH^*_{cont}(M)$.
Let $f$ be a Morse function and consider the complex of Novikov
chains
$$
CQ_{2n-k}(-\e f) = CM_{2n-k}(-\e f)\otimes \Lambda^\downarrow (=
CF_k(\e f)). \tag A.1
$$
On non-exact symplectic manifolds, this is typically {\it infinite
dimensional} as a $Q$-vector space. Therefore it is natural to put
some topology on it rather than to consider it just as an {\it
algebraic} vector space. For this purpose, we recall the
definition of the level $\lambda(\alpha)=\lambda_{\e f}(\alpha)$
of an element
$$
\alpha = \sum_A \alpha_A q^A, \quad \alpha_A \in CM_*(-\e f):
$$
$$
\align
\lambda(\alpha) & = \max\{\AA_{\e f}(\alpha_A q^A) \mid \alpha_A \neq 0 \} \\
& = \max\{\lambda^{Morse}_{- \e f}(\alpha_A) - \omega(A) \}.
\endalign
$$
As we saw before, this level function satisfies the inequality
$$
\lambda(\alpha + \beta) \leq \max\{\lambda(\alpha),
\lambda(\beta)\} \tag A.2
$$
and provides a natural filtration on $CQ_{2n-k}(-\e f)$, which
defines a {\it Non-Archimedan topology .} We refer to [Br] for a
nice exposition to the Non-Archimedean topology and geometry.

\proclaim{Definition \& Proposition A.1}  For each degree $*$,
consider the collection
$$
\BB = \bigcup_{\alpha \in CQ_*(-\e f), R \in \R} \{ U(\alpha, R)
\subset CQ_*(-\e f)  \}
$$
of the subsets  $U(\alpha, R)$ that is defined by
$$
U(\alpha,R) = \{ \beta \in CQ_*(-\e f) \mid \lambda(\beta -
\alpha) < R \}.
$$
Then $\BB$ satisfies the properties of a basis of topology.  We
equip $CQ_*(-\e f)$ with the topology generated by the basis
$\BB$.
\endproclaim
\demo{Proof} We need to prove that for any given $U(\alpha_1,
R_1)$ and $U(\alpha_2, R_2)$ with $U(\alpha_1, R_1) \cap
U(\alpha_2, R_2)\neq \emptyset$ and for any $\alpha \in
U(\alpha_1, R_1) \cap U(\alpha_2, R_2)$, there exists $R_3$ such
that
$$
U(\alpha,R_3) \subset U(\alpha_1, R_1) \cap U(\alpha_2, R_2). \tag
A.3
$$
Let $\beta \in U(\alpha_1, R_1) \cap U(\alpha_2, R_2)$. Then
$\beta$ satisfies
$$
\lambda(\beta - \alpha_i) < R_i, \quad i = 1, 2 \tag A.4
$$
Suppose $\gamma \in U(\beta, R)$ where $R$ is to be determined.
Then we derive from (A.2)
$$
\lambda(\gamma - \alpha_i) \leq \max\{ \lambda(\gamma - \beta),
\lambda(\beta - \alpha_i) \} = \max\{R, R_i\} \tag A.5
$$
Therefore if we choose $R \leq \min\{ R_1, R_2\}$, then we will
have
$$
U(\beta, R) \subset U(\alpha_1, R_1) \cap U(\alpha_2, R_2)
$$
which finishes the proof of the fact that $\BB$ really defines a
basis of topology. \qed\enddemo By the Non-Archimedean triangle
inequality (A.2), it follows that the basis element $U(\alpha,R)$
is nothing but  the affine subspace
$$
U(\alpha, R) = CQ_*^R(- \e f) + \alpha = CF_{2n -*}^R(\e f) +
\alpha
$$
where $CF_*^R$ is defined as in section 4.

The following is an easy consequence of the definition of the
boundary operator.

\proclaim{Lemma A.2} The  boundary operator
$$
\part_{\e f}= \part^{Morse}_{-\e f}\otimes \Lambda:
CQ_{2n-k}(-\e f) \to CQ_{2n-k-1}(-\e f)
$$
is continuous with respect to this topology.
\endproclaim
\demo{Proof} Let $U(\alpha,R)$ be a basis element and consider the
preimage
$$
(\part_{\e f})^{-1}(U(\alpha,R)).
$$
Suppose $\beta \in (\part_{\e f})^{-1}(U(\alpha,R))$, i.e.,
$\part_{\e f}(\beta) \in U(\alpha,R)$ and so
$$
\lambda(\part_{\e f}(\beta) - \alpha) < R. \tag A.6
$$
Recall that
$$
\lambda(\part_{\e f}(\delta))  \leq  \lambda(\delta) \tag A.7
$$
for any Novikov Floer chain $\delta$. Now we consider the basis
element $U(\beta, R)$. Then if $\gamma \in U(\beta,R)$, we have
$$
\aligned \lambda(\part_{\e f}(\gamma) - \alpha) & \leq
\max\{\lambda(\part_{\e f}(\gamma - \beta) ), \lambda(\part_{\e f}(\beta)- \alpha)\}\\
& \leq \max\{\lambda(\gamma - \beta), \lambda(\part_{\e f}(\beta)- \alpha)\}\\
& < \max\{R, R \} = R
\endaligned
\tag A.8
$$
where the second inequality comes from (A.7). This finishes the
proof of $\part_{\e f}(U(\beta, R)) \subset U(\alpha,R)$ i.e.,
$U(\beta, R) \subset (\part_{\e f})^{-1}(U(\alpha,R)$ for any
$\beta \in U(\alpha,R)$. Hence the proof. \qed\enddemo

Now we define

\definition{Definition A.3} A linear functional
$$
\mu: CQ_{2n-k}(-\e f) \to \Q
$$
is called continuous (or bounded) if it is continuous with respect
to the
 topology induced by the above filtration. We denote by
$CQ^\ell_{cont}(-\e f)$ the set of continuous linear functionals
on $CQ_{2n-k}(-\e f)$.
\enddefinition

The following is easy to see from the definition of  Novikov
chains. \proclaim{Lemma A.4} A linear functional $\mu$ is
continuous if and only if there exists $\lambda_\mu \in \R$ such
that
$$
\mu(\alpha_A q^A) = 0 \tag A.9
$$
for all $A$ with $-\omega(A) \leq \lambda_\mu$.
\endproclaim
\demo{Proof} The sufficiency  part of the proof is easy and so we
will focus on the necessary condition. We will prove this by
contradiction. Suppose that $\mu:CQ_{2n-k}(-\e f) \to \Q$ is a
continuous linear functional, but there exists a sequence of $A_j$
with
$$
-\omega(A_j) \to-\infty, \quad i.e.,  \, \, \omega(A_j) \to +
\infty \tag A.10
$$
and $\alpha_j \in CM_*(-\e f)$ such that
$$
\mu(\alpha_jq^{A_j}) \neq 0. \tag A.11
$$
Now consider the sequence of  Novikov chains
$$
\beta_N= \sum_{j=1}^N \alpha_j q^{A_j}. \tag A.12
$$
It is easy to check from (A.10) that $\beta_N$ converges to the
Novikov chain
$$
\beta = \sum_{j=1}^\infty \alpha_j q^{A_j}
$$
in the given Non-Archimedean topology on $CQ_*(-\e f)$. In fact,
this convergence holds for the sequence
$$
\beta_{c, N} = \sum_{j=1}^N (c_j \alpha_j) q^{A_j} \tag A.13
$$
for any given sequence $c =\{c_j \in \Q\}_{1 \leq j < \infty}$. We
choose $c_j$'s so that
$$
c_j = \frac{1}{\mu(\alpha_j q^{A_j}) }
$$
which is well-defined by (A.11).  However we  then have
$$
\mu(\beta_{c,N+1}) - \mu(\beta_{c,N})  =
\mu(c_{N+1}\alpha_{N+1}q^{A_{N+1}}) =
c_{N+1}\mu(\alpha_{N+1}q^{A_{N+1}}) =  1
$$
for all $N$. This proves that $\mu$ cannot be continuous, a
contradiction. This finishes the proof. \qed\enddemo

It then follows that
$$
\part_Q^* = \part_{-\e f}^*: (CQ_\ell(-\e f))^* \to (CQ_{\ell+1}(-\e
f))^*
$$
maps continuous linear functionals to continuous ones and so
defines the canonical complex
$$
(CQ^*_{cont}(-\e f), \part_Q^*)
$$
and hence defines the homology
$$
QH^\ell_{cont}(M): = H^\ell(CQ^*_{cont}(-\e f), \part_Q^*)).
$$
We recall the canonical embedding
$$
\sigma: CQ^\ell(-\e f)= CM_{2n-\ell} {\e f}\otimes
\Lambda^\uparrow \hookrightarrow CQ^\ell_{cont}(-\e f); a \mapsto
\langle a, \cdot \rangle \tag A.14
$$
mentioned in Remark 5.1. We have the following proposition which
is straightforward to prove. We refer to the proof of [Proposition
2.2, Oh2] for the details.

\proclaim{Proposition A.5} The map $\sigma$ in (A.14) is a chain
map from $(CQ^\ell(-\e f), \delta^Q)$ to $(CQ^\ell_{cont}(-\e f),
\part_Q^*)$. In particular we have a natural degree preserving
homomorphism
$$
\sigma: QH^*(M) \cong HQ^*(-\e f) \to HQ^*_{cont}(-\e f) \cong
QH^*_{cont}(M). \tag A.15
$$
\endproclaim
Now we can define the notion of continuous Floer cohomology
$HF^*_{cont}(H)$ for any given Hamiltonian in a similar way. Then
the co-chain map
$$
(h_\HH)^*: CF^k(H) \to CF^k(\e f)
$$
restricts to the co-chain map
$$
(h_\HH)^*: CF^k_{cont}(H) \to CF^k_{cont}(\e f).
$$
Once we have defined the continuous quantum cohomology and the
continuous Floer cohomology, it is straightforward to define the
spectral invariants for the continuous cohomology class in the
following way.

\definition{Definition A.6} Let $\mu \in QH^\ell_{cont}(M)$. Then we define
$$
\rho(H;\mu): = \inf \{\lambda \in \R \mid \mu \in \text{Im
}i_\lambda^* \} \tag A.16
$$
\enddefinition

Now it is straightforward to generalize all the axioms in Theorem
I to the continuous quantum cohomology class. The only non-obvious
axiom is the triangle inequality. But the proof will be a verbatim
modification of [Theorem II (5), Oh2] incorporating the argument
in the present paper that uses the Hamiltonian fibration and
pseudo-holomorphic sections. We leave the details to the
interested readers. We hope to investigate further properties of
the continuous quantum cohomology and its applications elsewhere.

\head {\bf References}
\endhead

\enddocument